\documentclass{amsart}

\usepackage{amssymb,amsmath,amsthm,color}

\usepackage[bookmarksopen=true,bookmarksnumbered=true, bookmarksopenlevel=2, colorlinks=true,linkcolor=mblue,
citecolor=mblue,urlcolor=mblue, pdfstartview=FitH]{hyperref}

\input xy
\xyoption{all}
\newcommand{\myar}{\ar@<0.5ex>@/^2ex/@{-}[rr]}
 \newcommand{\abcd}[4]{\left(
         \begin{smallmatrix}#1&#2\\#3&#4\end{smallmatrix}\right)}
\newcommand{\abcdm}[4]{\left(\begin{matrix}#1&#2\\
#3 & #4
\end{matrix}\right)}

\definecolor{mblue}{rgb}{0,0,.8}

\usepackage[all]{xy}

\newcommand{\N}{\mathbb N}
\newcommand{\Z}{\mathbb Z}
\newcommand{\Q}{\mathbb Q}
\newcommand{\F}{\mathbb F}

\newcommand{\PP}{\mathbb P}
\newcommand{\C}{\mathbb C}

\newtheorem{thm}{Theorem}
\newtheorem{lem}[thm]{Lemma}

\newtheorem{prop}[thm]{Proposition}
\newtheorem*{proposition}{Proposition}
\newtheorem{cor}[thm]{Corollary}

\newtheorem{ex1}[thm]{Example}
\newtheorem{defn}[thm]{Definition}
\newtheorem{rem}[thm]{Remark}

 \DeclareMathOperator{\SL}{SL}  
\DeclareMathOperator{\Gal}{Gal}   \DeclareMathOperator{\PSL}{PSL}
   \DeclareMathOperator{\sign}{sign}

\def\dash---{\thinspace---\hskip.16667em\relax}

\begin{document}

\title{Lifts of projective congruence groups}
\author{Ian Kiming, Matthias Sch\"{u}tt, Helena A.\ Verrill}

\address[Ian Kiming]{Department of Mathematical Sciences, University of Copenhagen, Universitetsparken 5, 2100 Copenhagen \O , Denmark}
\email{\href{mailto:kiming@math.ku.dk}{kiming@math.ku.dk}}
\address[Matthias Sch\"{u}tt]{Institut f\"ur Algebraische Geometrie, Leibniz Universit\"at Hannover, Welfengarten 1, 30167 Hannover, Germany}
\email{\href{mailto:schuett@math.uni-hannover.de}{schuett@math.uni-hannover.de}}
\address[Helena A.\ Verrill]{Department of Mathematics, Louisiana State University, Baton Rouge, Louisiana 70803-4918, USA}
\email{\href{mailto:verrill@math.lsu.edu}{verrill@math.lsu.edu}}

\begin{abstract}
We show that noncongruence subgroups of $\SL_2(\Z)$ projectively equivalent to congruence subgroups are ubiquitous. More precisely, they always exist if the congruence subgroup in question is a principal congruence subgroup $\Gamma(N)$ of level $N>2$, and they exist in many cases also for $\Gamma_0(N)$.

The motivation for asking this question is related to modular forms: projectively equivalent groups have the same spaces of cusp forms for all even weights whereas the spaces of cusp forms of odd weights are distinct in general. We make some initial observations on this phenomenon for weight $3$ via geometric considerations of the attached elliptic modular surfaces.

We also develop algorithms that construct all subgroups projectively equivalent to a given congruence subgroup and decide which of them are congruence. A crucial tool in this is the generalized level concept of Wohlfahrt.
\end{abstract}



\maketitle

\section{Introduction}

Suppose that $\Gamma_1$ and $\Gamma_2$ are subgroups of finite index of $\SL_2(\Z)$ that are projectively equivalent, i.e., have the same image in $\PSL_2(\Z)$. Thus $\langle \Gamma_1 , -1 \rangle = \langle \Gamma_2 , -1 \rangle$, and so the space of modular forms of given even weight is the same for the two groups $\Gamma_i$. But the spaces of forms of odd weights will in general be different for the two groups (see Section \ref{s:ex} and \ref{s:MF} for concrete examples and results).
\smallskip

The motivations behind the present paper are twofold: first the question whether such a situation can occur with $\Gamma_1$ a congruence subgroup but $\Gamma_2$ a noncongruence subgroup; and secondly, in such cases to study the attached spaces of cusp forms of odd weights further.
\smallskip

In the present paper we focus primarily on the first question and find that the answer is a resounding `yes'. More precisely, we give a complete answer to the question in case one of the groups is a principal congruence subgroup $\Gamma(N)$, as well as a partial answer for $\Gamma_0(N)$.
\smallskip

We employ the following terminology: For a subgroup $\Gamma$ of $\SL_2(\Z)$ denote by $\overline\Gamma$ the image of $\Gamma$ in $\PSL_2(\Z)$. By a {\it lift} of $\overline\Gamma$ we mean a subgroup of $\SL_2(\Z)$ that projects to $\overline\Gamma$ in $\PSL_2(\Z)$. A lift is called a {\it congruence lift} if it is a congruence subgroup.

As usual, if $\Gamma$ is a congruence subgroup then by the {\it level} of $\Gamma$ we understand the least $M$ such that $\Gamma(M)\subset\Gamma$.
\smallskip

Our main results are as follows.

\begin{thm}\label{Thm:gamma} Let $N\in\N$ and
$\Gamma(N)$ the principal congruence subgroup of $\SL_2(\Z)$ with
projective image $\overline{\Gamma(N)}\le \PSL_2(\Z)$.
\begin{enumerate}
\item
The number $n(N)$ of congruence lifts of $\overline{\Gamma(N)}$ is exactly
\[
 n(N) = \begin{cases}
         1 & \text{if } N=1,\\
3 & \text{if $N>1$ is odd},\\
5 & \text{if } N=2,\\
9 & \text{if $N>2$ is even}.
        \end{cases}
\]
\item
If $N>2$, then $\overline{\Gamma(N)}$ has noncongruence lifts.
\end{enumerate}
\end{thm}

\begin{thm}\label{Thm:gamma0} Let $p$ be a prime number and let $N\in\N$.
\smallskip

\noindent (i) If $4\nmid N$ and all odd prime divisors of $N$ are congruent to $1$ modulo $4$, then all lifts of $\overline{\Gamma_0(N)}$ are congruence.
\smallskip

\noindent (ii) If $p\equiv 3\pod{4}$ and $N=p^r$ for $r\in\N$, then there are precisely $3$ congruence lifts of $\overline{\Gamma_0(N)}$, namely $\Gamma_0(N)$, the subgroup of $\Gamma_0(N)$ consisting of those elements whose diagonal entries are squares modulo $p$ and one further subgroup (see Section \ref{Proof:gamma0}).

\noindent (iii) If $N$ is divisible by  $6, 9, 16, 20$ or by a prime $p>3$ congruent to $3$ modulo $4$, then the groups $\overline{\Gamma_0(N)}$ and $\overline{\Gamma_1(N)}$ have noncongruence lifts.
\end{thm}

We will give a detailed discussion of the group $\overline{\Gamma_0(N)}$ for $N=4,6,8,16,20$ in section \ref{examples}. It is shown that all lifts of $\overline{\Gamma_0(N)}$ are congruence if $N=4,8$, but that there are noncongruence lifts 
when $N=6,16,20$.
\smallskip

Thus, as far as the question of existence of noncongruence lifts of the group $\overline{\Gamma_0(N)}$ is concerned, Theorem \ref{Thm:gamma0} leaves undecided only the cases where $N$ is $3, 4$ or $8$ times an odd number whose prime divisors are all $\equiv 1 \pod{4}$.
\smallskip

A crucial tool is a variant of a result of Wohlfahrt \cite{wohlfahrt}. The difference between \cite[Theorem 2]{wohlfahrt} and the following proposition is that Wohlfahrt deals with subgroups of $\PSL_2(\Z)$ rather than $\SL_2(\Z)$. It was a minor surprise to us that his result did not just carry over literally, but that the level $N$ in the conclusion had to be replaced by $2N$. The notion of `general level' of a subgroup of finite index in $\SL_2(\Z)$ (or $\PSL_2(\Z)$) will be recalled in section \ref{level_concept}.

\begin{prop}\label{Prop:2N} Let $\Gamma$ be a subgroup of finite index in $\SL_2(\Z)$ of general level $N$.
\smallskip

If $\Gamma$ is congruence, then $\Gamma(2N)\le\Gamma$, and so $\Gamma$ has level $N$ or $2N$.
\end{prop}

The paper is organized as follows: in the first four sections, we establish the techniques that we will use to prove our main results.
First we elaborate in generality on lifts of subgroups from $\PSL_2(\Z)$ to $\SL_2(\Z)$.
Then we recall Wohlfahrt's generalized level concept and adjust it to our situation (Proposition \ref{Prop:level}).
In section \ref{s:alg}, we derive an algorithm that determines all lifts of a given subgroup of $\PSL_2(\Z)$ and decides which of the lifts are congruence (Proposition \ref{Prop:algorithm}).
This algorithm enables us to investigate some examples in detail in section \ref{s:ex}.

Sections \ref{Proof:gamma} and \ref{Proof:gamma0} continue to give the proofs of Theorems \ref{Thm:gamma} and \ref{Thm:gamma0}.
Next to the level concept and the examples, there is a crucial contribution by information on possible representations of the groups $\overline{\Gamma(N)}$ and $\overline{\Gamma_0(p)}$ ($p$ prime), due to Frasch \cite{frasch} and Rademacher \cite{rademacher}.
\smallskip

The paper concludes with additional observations: first, in section \ref{s:MF}, we notice that projective equivalence of subgroups of finite index in $\SL_2(\Z)$ not containing $-1$ is in fact equivalent to the condition that for infinitely many $k$ the groups have the same modular forms of weight $k$.

Secondly, in section \ref{squares} we draw some consequences of our main results for the groups generated by squares of elements in congruence subgroups.
We establish results when these groups are again congruence subgroups.

Finally, in section \ref{ss:EMS} we note some observations on spaces of cusp forms of weight $3$ for different lifts. Our arguments are geometric in nature, based on the elliptic modular surfaces attached to the lifts.

\section{Preliminaries on lifts}\label{prelim}

We shall first prove a simple but basic lifting lemma.

\begin{lem}\label{Lemma:lifts} Consider a subgroup $\overline\Gamma$ of $\PSL_2(\Z)$ of finite index. Suppose that we are given a presentation of $\overline\Gamma$ in terms of generators $\bar g_1,\ldots \bar g_s$ and relations $\bar R_1=1,\ldots ,\bar R_t=1$. The relations $\bar R_j$ have form:
$$
\bar R_j = \prod_{k=1}^{m_j} \bar h_{j,k}
$$
with each $\bar h_{j,k} \in \{ \bar g_1,\ldots \bar g_s \}$.
\smallskip

For $i=1,\ldots ,s$ and $j=1,\ldots ,t$ define the non-negative integer $\sigma_{i,j}$ to be the number of occurrences of $\bar g_i$ in the relation $\bar R_j$, i.e., the number of $k\in \{ 1,\ldots ,m_j \}$ such that $\bar h_{j,k} = \bar g_i$.
\medskip

\noindent (i). The group $\overline\Gamma$ has a lift $\Gamma$ in $\SL_2(\Z)$ such that $-1\not\in\Gamma$ if and only if the $\bar g_i$ have lifts $g_i\in \SL_2(\Z)$ such that:
$$
R_j = 1
$$
where for each $j=1,\ldots ,t$ the element $R_j$ is defined as the product $\prod_{k=1}^{m_j} h_{j,k}$ where $h_{j,k} := g_i$ if $\bar h_{j,k} = \bar g_i$.
\medskip

\noindent (ii). Suppose that $\overline\Gamma$ has a lift $\Gamma$ in $\SL_2(\Z)$ such that $-1\not\in\Gamma$, given by generators $g_1,\ldots g_s$ as in (i).

Then the lifts of $\overline\Gamma$ not containing $-1$ are parametrized by solutions
$$
(x_1,\ldots,x_s) \in \F_2^s
$$
to the linear system of equations
$$
(x_1,\ldots,x_s) ((\sigma_{i,j} \bmod{2}))_{i,j} = (0,\ldots,0)
$$
over $\F_2$. Here, a given solution $X=(x_1,\ldots,x_s)$ corresponds to the subgroup $\Gamma_X$ of $\SL_2(\Z)$ generated by:
$$
\{ g_i \mid ~x_i=0 \} \cup \{ -g_i \mid ~x_i=1 \} ~.
$$
\end{lem}

\begin{proof} Let $\alpha$ be the canonical homomorphism from $\SL_2(\Z)$ onto $\PSL_2(\Z)$.
\smallskip

We start by proving {\it (i)}. The necessity of the condition is clear: Suppose that $\overline\Gamma$ has a lift $\Gamma$ not containing $-1$. Let $g_i\in\Gamma$ be lifts of the $\bar g_i$. The products $R_j$ defined in the proposition then all map to $1$ in $\overline\Gamma$ under $\alpha$, and hence $R_j = \pm 1$ for all $j$. However, the $R_j$ are also elements of $\Gamma$, so we must have $R_j=1$ for all $j$.
\smallskip

Conversely, suppose that there are lifts $g_i$ such that $R_j=1$ for $j=1,\ldots ,t$. Then let $\Gamma$ be the subgroup of $\SL_2(\Z)$ generated by the $g_i$.

Since we have the relations $R_j=1$ among these generators, $\Gamma$ must be -- as an abstract group -- a quotient of $\overline\Gamma$. Let $\beta : ~ \overline\Gamma \rightarrow \Gamma$ be the corresponding homomorphism. It is given by $\beta(\bar g_i) = g_i$.

On the other hand, $\Gamma$ also surjects to $\overline\Gamma$ via $\alpha$. We have $\alpha(g_i) = \bar g_i$.

Now, $(\alpha \circ \beta)(\bar g_i) = \bar g_i$ for each $i$ and this shows that $\alpha \circ \beta = \mbox{id}$. Since $\beta$ is a surjection onto $\Gamma$, we conclude that $\alpha$ is injective on $\Gamma$. Thus, $-1\not\in \Gamma$.
\smallskip

For the proof of part {\it (ii)} we first note that it is clear from the above that any lift of $\overline\Gamma$ has shape $\Gamma_X$ for some solution $X$ to the stated system of equations over $\F_2$, and that the group $\Gamma_X$ is in fact a lift of $\overline\Gamma$ for any such solution $X$. To finish the proof, note that any such $\Gamma_X$ does not contain $-1$; this follows because $-1\not\in\Gamma$. For the same reason, we also see that the subgroups corresponding to two distinct solutions $X$ are in fact distinct subgroups.
\end{proof}

We will also need the following observation concerning subgroups of projective groups with noncongruence lifts.

\begin{lem}
\label{Lem:sub}
Let $\overline\Gamma$ be a subgroup of $\PSL_2(\Z)$ with a noncongruence lift to $\SL_2(\Z)$.
Then any subgroup $\overline{G}\le \overline\Gamma$ has a noncongruence lift.
\end{lem}

\begin{proof}
Let $\Gamma$ denote a noncongruence lift of $\overline\Gamma$.
Define $G$ as the pre-image of $\overline{G}$ under the projection $\Gamma\to\overline\Gamma$.
By construction, $G$ is a lift of $\overline{G}$ to $\SL_2(\Z)$.
As a subgroup of the noncongruence subgroup $\Gamma$, $G$ cannot be congruence.
\end{proof}

\section{Level concept}\label{level_concept}

In this section, we will derive the following restrictions on the level of congruence lifts:

\begin{prop}\label{Prop:level} Let $N\in\N$ and let $\overline\Gamma$ be a subgroup of $\PSL_2(\Z)$ with
$$
\overline{\Gamma(N)} \le \overline\Gamma \le \overline{\Gamma_0(N)} ~.
$$

Then any congruence lift of $\overline\Gamma$ has level $N$ or $2N$.
\end{prop}

In section \ref{examples}, we will show by means of an example that the proposition cannot be improved: both possibilities $N$ and $2N$ for the level do in fact occur.
\smallskip

To prepare the proof of Proposition \ref{Prop:level} we will first recall the generalized notion of level for arbitrary subgroups of finite index in $\SL_2(\Z)$, as introduced by Wohlfahrt, cf.\ \cite{wohlfahrt}.
\smallskip

Let $\Gamma$ be a subgroup of finite index in $\SL_2(\Z)$. The general level $N$ of $\Gamma$ is defined as the least common multiple of the cusp widths. Recall the width of a cusp $c$ of $\Gamma$ is the least possible $n\in\N$ such that $\pm g T^n g^{-1} \in \Gamma$ where $g\in\SL_2(\Z)$ is such that $g\infty = c$, and where $T$ here as well as throughout the paper denotes the usual translation matrix
$$
T:=\begin{pmatrix}1 & 1\\0 & 1\end{pmatrix} .
$$

Notice that the general level of $\Gamma$ only depends on $\overline\Gamma$. Thus, the general level of any lift of $\overline\Gamma(N)$ is $N$. It is easy to see that if $\Gamma_1$ and $\Gamma_2$ are two subgroups of finite index in $\SL_2(\Z)$, and if $\Gamma_2 \le \Gamma_1$, then the general level of $\Gamma_1$ is a divisor of the general level of $\Gamma_2$.
\smallskip

\begin{lem}\label{lem1} Let $\Gamma$ be a subgroup of finite index in $\SL_2(\Z)$. Let $N$ denote the general level of $\Gamma$.
\smallskip

\noindent (a). If $\nu\in\N$ is such that $\Gamma(\nu)\le\Gamma$, then $N\mid \nu$.
\medskip

\noindent (b). Let $N\mid \nu$. Then for any $g\in \SL_2(\Z)$, we have $gT^{2\nu}g^{-1}\in\Gamma$.
\end{lem}

\begin{proof} {\it (a)} By the inclusion $\Gamma(\nu)\le\Gamma$, the width with respect to $\Gamma$ of any cusp divides the width with respect to $\Gamma(\nu)$ -- which is $\nu$. The claim follows.
\medskip

\noindent {\it (b)} Suppose that $N\mid \nu$ and consider any $g\in\SL_2(\Z)$. The stabilizer in $\overline\Gamma$ of the cusp $g^{-1}\infty$ is $\{ \bar{A}^m \mid ~ m\in\Z \}$ for a certain matrix $A\in\Gamma$.

Then $g^{-1}Ag$ stabilizes $\infty$, and by definition of the width $b$ of the cusp $g^{-1}\infty$, we have $g^{-1}Ag = \pm T^b$. By definition of the general level $N$ we have $b\mid N$, hence $b\mid \nu$, and so $\pm gT^{\nu}g^{-1} = A^{\nu / b} \in \Gamma$. Thus $gT^{2\nu}g^{-1} \in \Gamma$ as desired.
\end{proof}

Let us recall the statement of Proposition \ref{Prop:2N}. The proposition is a slight variant of a result of Wohlfahrt, cf.\ Theorem 2 of \cite{wohlfahrt} and the succeeding remark.

\begin{proposition} Let $\Gamma$ be a subgroup of finite index in $\SL_2(\Z)$ of general level $N$.
\smallskip

If $\Gamma$ is congruence, then $\Gamma(2N)\le\Gamma$, and so $\Gamma$ has level $N$ or $2N$.
\end{proposition}

\begin{proof} Suppose that $\nu\in\N$ is such that $\Gamma(\nu)\le\Gamma$ and let $A=\begin{pmatrix} a & b\\c & d \end{pmatrix}\in\Gamma(2N)$ be arbitrary. We must show that $A\in\Gamma$. Note that it will be enough to show $BAC\in\Gamma$ with matrices $B,C\in \Gamma \cap \Gamma(2N)$, so we can modify $A$ in this way whenever convenient.

We now mimic some computations by Wohlfahrt, cf.\ \cite{wohlfahrt}, proof of Theorem 2 (and in the process also correct a small typo in his argument).
\smallskip

First we claim that we can modify $A$ in the above manner so as to obtain $(d,\nu)=1$. If $d=\pm 1$ this is already the case, so assume $d\not= \pm 1$. Then $c\not=0$, but $(c,d)=1$. Since $A\in\Gamma(2N)$ we have $(d,2N)=1$ so that $(d,2Nc)=1$. As now $2Nc\not= 0$, we deduce by Dirichlet's theorem on primes in arithmetic progressions that there is $m\in\N$ such that $d+m\cdot2Nc$ is a prime larger than $\nu$. In particular, this implies $(d+2Nmc,\nu) = 1$. By Lemma \ref{lem1} {\it (b)} we have $T^{2N} \in\Gamma \cap \Gamma(2N)$, and so we may replace $A$ by $AT^{2Nm}$ and thus assume that $(d,\nu)=1$ (and still $A\in\Gamma(2N)$).

Since now $(d,\nu)=1$ and $N\mid \nu$ by Lemma \ref{lem1} {\it (a)}, we have $(Nd,\nu) = N$. Consequently, $(2Nd,\nu)$ is a divisor of $2N$ which in turn divides $b$ as $A\in\Gamma(2N)$. Hence the congruence $b+n\cdot 2Nd \equiv 0 \pod{\nu}$ is solvable for $n$. Replacing $A$ by $T^{2Nn}A$ does not change $d$, but changes $b$ to $b+2Nnd$. So we may additionally assume $b\equiv 0 \pod{\nu}$.

Now, again by Lemma \ref{lem1} {\it (b)} the matrix
$$
\begin{pmatrix} 1 & 0 \\ 2Nm & 1 \end{pmatrix} = \begin{pmatrix} 0 & -1 \\ 1 & 0 \end{pmatrix} T^{-2Nm} \begin{pmatrix} 0 & 1 \\ -1 & 0 \end{pmatrix}
$$
is in $\Gamma \cap \Gamma(2N)$ for any $m\in\Z$. With an argument similar to the above we see that if we multiply $A$  on the left by this matrix for a suitable $m$, we may assume $c\equiv 0 \pod{\nu}$. This multiplication leaves $b$ unchanged, and the condition $(d,\nu)=1$ is preserved since $b\equiv 0 \pod{\nu}$.

We now have $b\equiv c\equiv 0 \pod{\nu}$ and $(d,\nu)=1$. Then necessarily $ad\equiv 1\pod{\nu}$.
\smallskip

Now consider the matrix:
\begin{eqnarray}\label{eq:M}
M = (T^{1-d})' \underbrace{\begin{pmatrix} a & a-1\\1-a & 2-a \end{pmatrix}}_{L} T^{d-1} = \begin{pmatrix} a & ad-1 \\ 1-ad & d(2-ad)\end{pmatrix}
\end{eqnarray}
where we have denoted by $D'$ the transpose of $D$ for any matrix $D$.

We find $A\equiv M\mod \nu$, so that $M^{-1}A\in\Gamma(\nu) \le \Gamma$.
To deduce the claim $A\in\Gamma$, it thus suffices to show that $M\in\Gamma$.
By Lemma \ref{lem1} {\it (b)}, $(T^{1-d})'$ and $T^{d-1}$ are both in $\Gamma$, since $d\equiv 1\mod 2N$ and $(T^{1-d})'$ is a conjugate of $T^{d-1}$.
Thus it suffices to verify that $L\in\Gamma$.
Here we factor
$$
L= T' T^{a-1} T'^{-1} ~.
$$

Since $a\equiv 1\pod{2N}$, Lemma \ref{lem1} {\it (b)} shows that $L\in\Gamma$. Hence we deduce that $M\in\Gamma$ and therefore $A\in\Gamma$ as claimed.
\medskip

Finally note that the level of $\Gamma$ is now seen to be $N$ or $2N$ by Lemma \ref{lem1} {\it (a)}.
\end{proof}

\begin{rem} In \cite{wohlfahrt} Wohlfahrt defines a congruence subgroup of $\SL_2(\Z)$ as a subgroup containing $\langle \pm 1, \Gamma(m) \rangle$ for some $m\in\N$.
His Theorem 2 states that if $\Gamma$ is a congruence subgroup, then $\Gamma \ge \langle \pm 1, \Gamma(N) \rangle$ with $N$ the general level of $\Gamma$.
Wohlfahrt's proof essentially coincides with the above proof of Proposition \ref{Prop:2N} with `$2N$ replaced by $N$'.
The difference between these cases lies in the fact that if $-1\in\Gamma$, then the conclusion of part {\it (b)} of Lemma \ref{lem1} can be improved to $gT^{\nu}g^{-1}\in\Gamma$ as an inspection of the proof immediately reveals.
There was a typo in Wohlfahrt's proof that we alluded to above: the matrix $\left( \begin{smallmatrix} a & ad-1 \\ 1-ad & d(2-ad)\end{smallmatrix} \right)$ does not equal
$$
(T^{d-1})' \begin{pmatrix} a & a-1\\1-a & 2-a \end{pmatrix} T^{d-1},
$$
but factors correctly as in (\ref{eq:M}).
\end{rem}

\begin{rem} The results in this section are also related to results of Larcher, cf.\ \cite{larcher}. However, he works with a different notion of `congruence subgroup', namely $-1$ is always assumed to be in the group in contrast to our situation.
\end{rem}

\begin{proof}[Proof of Proposition \ref{Prop:level}] Since all lifts of $\overline\Gamma$ have the same general level it will suffice by Proposition \ref{Prop:2N} to display a lift of general level $N$.

But $\overline\Gamma$ has a lift $\Gamma$ that sits between $\Gamma(N)$ and $\Gamma_0(N)$. Thus the general level of $\Gamma$ is a divisor of the general level of $\Gamma(N)$ and is a multiple of the general level of $\Gamma_0(N)$. The claim now follows since $\Gamma(N)$ and $\Gamma_0(N)$ both have general level $N$: all cusps of $\Gamma(N)$ have width $N$, and the cusp widths of $\Gamma_0(N)$ are all divisors of $N$, but the cusp $0$ has in fact width $N$.
\end{proof}

\section{Algorithm}
\label{s:alg}

In the next section, we will investigate two examples that will be used in the proof of Theorem \ref{Thm:gamma0}.
The basis for those examples are the following proposition and its corollary which we believe to be of independent interest.

\begin{prop}\label{Prop:algorithm}

Suppose that a subgroup $\overline\Gamma \le \PSL_2(\Z)$ of finite
index is given by a Farey symbol, or that there is a method of
determining whether an element of $\PSL_2(\Z)$ is in $\overline\Gamma$.

Then there is an algorithm that determines all lifts of $\overline\Gamma$ to $\SL_2(\Z)$ and decides which of the lifts are congruence.
\end{prop}

\begin{cor}
There is an algorithm that determines all congruence and noncongruence subgroups of $\SL_2(\Z)$ that are projectively equivalent to a given congruence subgroup.
\end{cor}

\subsection{Farey symbols}
\label{subsec:Farey}

Our algorithm is based on the use of Farey symbols, described in
\cite{kulkarni}.  We recall the definition, with a slight modification,
which is necessary because we work with subgroups of $\SL_2(\Z)$,
whereas Farey symbols as defined in \cite{kulkarni} are associated with
subgroups of $\PSL_2(\Z)$.

Recall Kulkarni's definition \cite{kulkarni} of a Farey symbol:
\begin{defn}
\label{def:FS}
A Farey symbol consists of a sequence
of cusps (elements of $\Q\cup\{\infty\}$),
of length $n+1$ for some positive integer
$n\ge 1$,
starting and ending with infinity,
(the starting term can be considered as $-\infty=\frac{-1}{0}$
and the last term
as $+\infty=\frac{1}{0}$)
 together with a sequence of
labels $l_1,l_2,\dots,l_n$ of length $n$.
The sequence of cusps,
$\{-\infty,\frac{a_2}{b_2},\frac{a_3}{b_3},\cdots,\frac{a_{n-1}}{b_{n-1}},
\infty\}$ satisfies $a_{i+1}b_i - a_ib_{i+1}=1$ for $i=1,\dots,n-2$,
and $b_1=b_{n-1}=1$.
Note that we take $a_1=-1,b_1=0$ and $a_{n+1}=1,b_{n+1}=0$.
A label is either a positive integer, or a
symbol $\circ$ or $\bullet$. For any integer occurring as a label, there
are exactly two labels having this value.
Such a symbol is written thus:
$$
\{
\xymatrix@R=0.5ex@C=0.5ex{
-\infty     \myar^{l_1}      &&
\frac{a_{2}}{b_{2}}           \myar^{l_2}   &&
\frac{a_{3}}{b_{3}} \myar^{l_3}       &&
{\phantom{\frac{}{}}}
\cdots
{\phantom{\frac{}{}}}
\myar^{l_{n-2}}       &&
\frac{a_{n-1}}{b_{n-1}}           \myar^{l_{n-1}}      &&
\frac{a_{n}}{b_{n}}           \myar^{l_n}        &&
\infty
}
\}
$$
\end{defn}

\subsection{Signed Farey symbols}

Farey symbols correspond to subgroups of $\PSL_2(\Z)$.
Since we are working with subgroups of $\SL_2(\Z)$ not containing
$-1$, we need a minor modification to Farey symbols.

\begin{defn}
A {\emph signed Farey symbol}
is defined exactly as a Farey symbol in Definition~\ref{def:FS}
above, except that the labels can be any nonzero integers, in pairs,
or the symbol $\bullet$.
The symbol $\circ$ is not used.
\end{defn}
Following Kulkarni, with a minor modification, the
subgroup of $\SL_2(\Z)$
corresponding to a signed Farey symbol
has generators described
as follows:
\begin{itemize}
\item
For any integer pair of integer labels $l_i=l_j$, with
$i<j$, we have
a generator
\begin{equation}
\label{eqn:farey_gij}
g_{ij}=\sign(l_i)
\begin{pmatrix}
a_jb_i + a_{j+1}b_{i+1} & -(a_{i+1}a_{j+1} + a_ia_j)\\
b_ib_j + b_{i+1}b_{j+1} & -(a_ib_j + a_{i+1}b_{j+1})
\end{pmatrix}.
\end{equation}
Note that it is important to fix the order $i<j$, since changing the order
changes the result from a matrix $m$ to $-m^{-1}$.

\item
For any label $l_j=\bullet$ we have an order $3$ generator
\begin{equation}
\label{eqn:farey_gj}
g_j=
\epsilon_j
\begin{pmatrix}
a_jb_j + a_jb_{j+1} + a_{j+1}b_{j+1} & -(a_{j+1}^2 + a_ja_{j+1} + a_j^2)\\
b_j^2 + b_{j}b_{j+1} + b_{j+1}^2 & -(a_jb_j + a_{j+1}b_j + a_{j+1}b_{j+1})
\end{pmatrix}
\end{equation}
where $\epsilon_j$ is chosen so that the matrix has order $3$
in $\SL_2(\Z)$ (rather than $6$).
The only relations
are those stating that generators corresponding
to $\bullet$ have order $3$.
\end{itemize}

Note that for the original Farey symbols,
there are additional generators corresponding to
labels $l_i=\circ$, namely
\begin{equation}
\label{eqn:farey_gi}
g_i=
\begin{pmatrix}
a_ib_i + a_{i+1}b_{i+1} & -(a_{i}^2 + a_{i+1}^2)\\
b_i^2 + b_{i+1}^2  & -(a_ib_i + a_{i+1}b_{i+1})
\end{pmatrix}.
\end{equation}
and additional relations, saying that this generator
has order $4$.  However, such groups necessarily contain
$-1$, and so there is no need to consider signs.

\begin{prop}
Any subgroup of finite index in $\SL_2(\Z)$ not containing $-1$ can be defined by
a signed Farey symbol, i.e., can be generated by the matrices
determined by such a symbol.
\end{prop}
\begin{proof}
This follows immediately from \cite[(5.4)]{kulkarni},
Lemma~\ref{Lemma:lifts}, and the fact that edges of Farey symbols
labeled $\circ$ correspond to matrices in $\SL_2(\Z)$ with
square $-1$, and so do not need to be included in this case.
\end{proof}

If a subgroup $\Gamma$ of finite index in
 $\SL_2(\Z)$ is given by specifying some rule for
determining whether or not an element of $\SL_2(\Z)$ is in the group,
then the Farey Symbol for $\Gamma$ is determind as described by
\cite{CLLT}, and then signs are determined by checking whether or not
each of the above matrices are in $\Gamma$.  If for any $i, j$ both
$g_{ij}$ and $-g_{ij}$ are in $\Gamma$, then also $-1$ is in $\Gamma$,
and so we do not specify any signs for the group.
(Similarly for $g_i$ or $g_j$.)

\subsection{Proof of Proposition \ref{Prop:algorithm}}


Suppose $\overline\Gamma$ is a subgroup of finite index in
$\PSL_2(\Z)$.  We want to determine all congruence lifts of $\Gamma$.

Lemma~\ref{Lemma:lifts} gives an algorithm to determine all lifts of $\Gamma$.

Suppose $\overline\Gamma$ can be given by a Farey symbol,
as described in the previous section,
using the method of Kulkarni \cite{kulkarni}.  For example, this holds true
if we can determine whether an element of $\PSL_2(\Z)$ is in $\overline\Gamma$.

In \cite{LLT}, Lang, Lim and Tan give an algorithm for determining
whether $\overline\Gamma$ is a projective congruence subgroup.
In the case that $\overline\Gamma$ is not a congruence subgroup, then
all lifts of $\Gamma$ are also noncongruence, and no more remains to be
done.

Suppose that $\overline\Gamma$ is determined to be a congruence
subgroup of level $N$ by the algorithm in \cite{LLT}.
The algorithm uses a Farey symbol for $\overline\Gamma$, and associated
generators $\overline g_1,\dots,\overline g_s$, and a set of generators
$h_1,\dots,h_k$ for $\Gamma(M)$, where $M$ is the general level $N$ of
$\overline\Gamma$, and if possible, determines a word for each $h_i$ in
terms of the $g_j$.  The same algorithm can also be used when $M=2N$.

Note that generators and relations of any principal congruence subgroup $\Gamma(M)$ can be obtained algorithmically by using the method of either
\cite{chuman} or \cite{CLLT} applied to the projective
image $\overline{\Gamma(M)}$, followed by
a simple check to determine sign choices of lifts of generators, and
application of
Lemma~\ref{Lemma:lifts} to determine relations.

Suppose that using this method, and once we have determined that
$\bar{h_i} = \prod \bar h_{i,j}$ for some $\bar h_{i,j}\in\{\bar g_1,
\dots,\bar g_s \}$,
then by checking signs, we have
$$h_i = (-1)^{\epsilon_i}\prod h_{i,j},$$
for $i=1,\dots,k$ for some $\epsilon_i\in\{0,1\}$,
where each $h_{i,j}\in\{g_1,\dots,g_s\}$.
Let $\sigma_{i,j}$ be the number of times $g_i$ occurs in this
product, modulo $2$,
i.e.,
\begin{equation}
\label{eqn:in_alg_sigmaij}
\sigma_{i,j}=\#\{i \;|\; h_{i,j}=g_j\} \mod 2.
\end{equation}
Then a lift of $\Gamma$ having
generators $(-1)^{\delta_i} g_i$ for $i=1,\dots s$
contains $\Gamma(M)$ if and only if
the vector
$(\delta_i)$ is a solution to the $\mod 2$ linear system of equations
\begin{equation}
\label{eqn:liftingcondition}
(\sigma_{i,j})(\delta_1,\cdots,\delta_s)=(\epsilon_i).
\end{equation}
By Proposition~\ref{Prop:level},
we need to test whether this is the case for $M=N$ and $M=2N$.
If $\Gamma(N)$ is contained in $\Gamma$, it is congruence of level $N$.
If $\Gamma(N)$ is not contained in $\Gamma$, but $\Gamma(2N)$ is,
then $\Gamma$ is congruence of level $2N$.  Otherwise
$\Gamma$ is not a congruence subgroup.
\qed

\begin{rem}
\label{rem:notesonalg}
In practice, there will be many generators of $\Gamma(M)$.
Let $V_{\Gamma,M}$ be the space corresponding to writing generators for
$\Gamma(M)$ as words in terms of the generators
for $\Gamma$. So $V_{\Gamma,M}$ is spanned by the vectors in $\F_2^s$,
\begin{equation}
\label{eqn:v_i_in_alg}
v_i=(\sigma_{i,1},\dots,\sigma_{i,s})
\end{equation}
 for $i=1,\dots,k$, where $\sigma_{i,j}$ are as in (\ref{eqn:in_alg_sigmaij}).
We only need take a set of $v_i$ spanning this
space in the system (\ref{eqn:liftingcondition}).

The signs simplify when
we have a known congruence lift of level $M$ not containing $-1$,
for example, if $\Gamma$ is already such a group.
Notably this happens for
$\Gamma_1(N)$, or
if we are considering lifts of $\overline\Gamma
=\overline{\Gamma_0(p)}$ for $p$ prime,
we can take the lift
$\Gamma$ consisting of matrices in $\Gamma_0(p)$ with
diagonal entries squares modulo $p$.  In this situation, we can take
all $\epsilon_i=0$.
\end{rem}

\begin{rem} Suppose that $\overline\Gamma$ has index $\mu$ in $\PSL_2(\Z)$, has $\nu_2$ elliptic points of order $2$, and has $\nu_3$ elliptic points of order $3$.

By Euler's formula applied to the tiling of the fundamental domain for $\overline\Gamma$ by images of the fundamental domain for $\PSL_2(\Z)$ corresponding to a Farey symbol for $\overline\Gamma$, the minimal number of generators of $\overline\Gamma$ is
$$
\delta := \frac{\mu}{6} + 1 + \frac{\nu_2}{2} + \frac{\nu_3}{3} ~.
$$

It follows that if $\nu_2=0$ then there are $2^{\delta - \nu_3}$ lifts of $\overline\Gamma$ not containing $-1$.
\end{rem}

\section{Examples}
\label{examples}
\label{s:ex}

The examples in this section were computed using Pari \cite{pari}, Magma \cite{magma} and GAP \cite{gap}.
The Magma program has built in functions for determining Farey symbols of congruence subgroups, as described in \cite{kulkarni}.
A GAP package \cite{DoJeKoVe} was also used for computing Farey symbols,
and working with subgroups of $\PSL_2(\Z)$ given only by their
Farey symbol.  The algorithm of \cite{LLT} for determining
whether some element of $\PSL_2(\Z)$ is in some
group given by a Farey symbol was also implemented in GAP.

\smallskip

In some of these examples, we also compute the dimensions of
some spaces of cusp forms.
The dimensions follow from Shimura's formula which we recall in Section~\ref{s:MF}.
They are given in terms of the number of regular and irregular
cusps, $\nu_\infty^+$ and $\nu_\infty^-$
since the subgroups in consideration have no torsion.

\begin{ex1}
\label{ex:4}
Let $\Gamma=\Gamma_1(4)$.
A signed
Farey symbol with corresponding fundamental domain $\mathcal F$
and generators $T, A$ is given in the following diagram:
$$
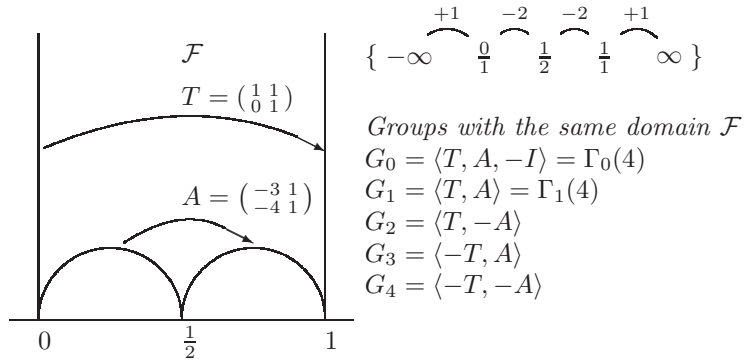
\begin{array}{ll}
\setlength\unitlength{1.5in}
\begin{picture}(1.1,1)(-0.1,-0.1)
\put(0,0){\line(0,1){1}}
\put(1,0){\line(0,1){1}}
\put(-0.1,0){\line(1,0){1.2}}
\qbezier(0,0)(0,0.103553)(0.07323,0.17677)
\qbezier(0.07323,0.17677)(0.14644,0.25)(0.25,0.25)
\qbezier(0.5,0)(0.5,0.103553)(0.42677,0.17677)
\qbezier(0.42677,0.17677)(0.353553,0.25)(0.25,0.25)
\qbezier(0.5,0)(0.5,0.103553)(0.57323,0.17677)
\qbezier(0.57323,0.17677)(0.64644,0.25)(0.75,0.25)
\qbezier(1,0)(1,0.103553)(0.92677,0.17677)
\qbezier(0.92677,0.17677)(0.853553,0.25)(0.75,0.25)
\qbezier(0.3,0.27)(0.5,0.4)(0.65,0.32)
\qbezier(0.02,0.6)(0.5,0.8)(0.90,0.64)
\put(0.90,0.64){\vector(2,-1){0.1}}
\put(0.65,0.32){\vector(2,-1){0.1}}
\put(0,-0.1){$0$}
\put(0.5,-0.1){$\frac{1}{2}$}
\put(1,-0.1){$1$}
\put(0.5,0.75){$T=\abcd{1}{1}{0}{1}$}
\put(0.5,0.4){$A=\abcd{-3}{1}{-4}{1}$}
\put(0.5,0.9){$\mathcal F$}
\end{picture}
&
\raisebox{1in}{$
\begin{array}{l}
\{
\xymatrix@R=0.5ex@C=0.5ex{
-\infty    \myar^{+1}      &&
\frac{0}{1}\myar^{-2}   &&
\frac{1}{2}\myar^{-2}      &&
\frac{1}{1}\myar^{+1}        &&
\infty
}
\}
\\
\\
\text{Groups with the same domain }\mathcal F\\
G_0=\langle T,A,-I\rangle=\Gamma_0(4)\\
G_1=\langle T,A\rangle=\Gamma_1(4)\\
G_2=\langle T,-A\rangle\\
G_3=\langle -T,A\rangle\\
G_4=\langle -T,-A\rangle\\
\end{array}
$
}
\end{array}
$$

We now apply the algorithm described in
Proposition~\ref{Prop:algorithm}. The generators for $\Gamma(4)$ in terms of $A$ and $T$ are as follows:
$$
\begin{array}{lllll}
\rule{0ex}{5ex}
    \abcdm{1}{4}{0}{1}  &=T^4  & \abcdm{5}{4}{-4}{-3}&=T^{-2}AT  \\
\rule{0ex}{5ex}
    \abcdm{1}{ 0}{4}{1}    &=A^{-1}T   &    \abcdm{5}{-4}{4}{ -3} &=TA^{-1}    \\
\rule{0ex}{5ex}
    \abcdm{-7}{12}{4}{-7}&=T^{-2}A^{-1}T^{-1}&\\
\end{array}
$$
see the proof of Lemma \ref{Lem:2^2} for more details on the generators of $\Gamma(4)$.

By the algorithm of Proposition~\ref{Prop:algorithm}, we see that level $4$ congruence lifts of $\overline{\Gamma_1(4)}$ correspond to solutions of
$
\begin{pmatrix}
1 & 1
\end{pmatrix}
\begin{pmatrix}
\delta_1 \\ \delta_2
\end{pmatrix}
=
\begin{pmatrix}
0 & 0\\
\end{pmatrix}
$,
so
$\langle T,A\rangle$
and
$\langle -T,-A\rangle$ are congruence subgroups of level $4$.

Using Magma \cite{magma} we compute that
$\Gamma(8)$ can be generated by $33$ matrices.  Writing these in terms of
$A$ and $T$, we find that each generator is given as a product of an
even number of $A$'s and $T$'s,
i.e., the matrix in (\ref{eqn:liftingcondition}) is the zero matrix,
and so the remaining lifts
$\langle -T,A\rangle$ and
$\langle T,-A\rangle$ must be
congruence subgroups of level $8$.

Thus all four groups $\langle \pm A,\pm T\rangle$ are
congruence subgroups.  Being congruence subgroups, they can also be
described by congruence conditions by considering the quotients
by $\Gamma(4)$ or $\Gamma(8)$.

The table below shows the regularity of the cusps of each lift of
$\overline\Gamma_0(4)$, which allows us to compute the
dimension of of $S_3(\Gamma)$ and $S_5(\Gamma)$ for each of these groups.
In this table $S=\abcd{0}{1}{-1}{0}$ and
$B=\abcd{1}{0}{2}{1}$.
$$
\begin{array}{|lllllll|}
\hline
\text{cusp}&\text{stabilizer}&\text{width}& \multicolumn{4}{c|}
{\text{regular in...?}}\\
&&&
G_1&
G_2&
G_3&
G_4\\
\infty & T & 1& Y & Y & n & n\\
0      & A^{-1}T= 
ST^{-4}S^{-1} &4 & Y & n & n & Y\\
\frac12    & A=-BT^{-1}B^{-1} & 1 & n & Y & n & Y\\
\multicolumn{3}{|r}{\dim S_3(\Gamma)}&0 & 0 & 1 & 0\\
\multicolumn{3}{|r}{\dim S_5(\Gamma)}&1 & 1 & 2 & 1\\
\multicolumn{3}{|r}{\text{level as congruence subgroup}}&4 &
8&8&4\\
\hline
\end{array}
$$


\end{ex1}

\begin{ex1}
\label{ex1:6}
The group
$\Gamma_1(6)$ has generators
$$T=\abcd{1}{1}{0}{1},\;
A=\abcd{-5}{1}{-6}{1},\;
B=\abcd{7}{-3}{12}{-5},$$
corresponding to a Farey symbol
$$
\{
\xymatrix@R=0.5ex@C=0.5ex{
-\infty    \myar^{+1}      &&
\frac{0}{1}\myar^{-2}   &&
\frac{1}{3}\myar^{+3}       &&
\frac{1}{2}\myar^{+3}      &&
\frac{2}{3}\myar^{-2}        &&
\frac{1}{1}\myar^{+1}        &&
\infty
}
\}
$$
%
There are no relations between these matrices, so
there are $8$ possible subgroups which are lifts of $\Gamma_1(6)$
not containing $-I$.  By Proposition~\ref{Prop:level}
if these are congruence, they have level $6$ or $12$.

Using Magma \cite{magma},
we find that $\Gamma(6)$ can be generated by the following 13 matrices.
The algorithm of \cite{LLT} expresses these generators in terms of
$g_1=A,g_2=B,g_3=T$.  In the table we also write the corresponding
vectors $v_i$, as  described in (\ref{eqn:v_i_in_alg})
of Proposition~\ref{Prop:algorithm}.
$$
\begin{array}{llll}
\abcd{-11}{  6}{ -24}{ 13  }&=B^{-2}\rule{0ex}{4ex} &(0,0,0)\\
\abcd{-17}{  6}{ -54}{ 19  }&=(A^{-1}B)^3  \rule{0ex}{4ex}&(1,1,0)\\
\abcd{-29}{ 18}{-108}{ 67  }&=A^{-1}BAB^{-1}\rule{0ex}{4ex}&(0,0,0)\\
\abcd{19}{ -12}{ 84}{ -53  }&=(BA^{-1}BA)^{-1}\rule{0ex}{4ex}&(0,0,0)\\
\abcd{-17}{ 12}{ -78}{ 55  }&=A^{-2}BA^{-1}\rule{0ex}{4ex}&(1,0,0)\\
\abcd{-41}{ 30}{-108}{ 79  }&=B^{-1}ABA^{-1}  \rule{0ex}{4ex}&(0,0,0)\\
\abcd{-23}{ 18}{ -78}{ 61  }&=A^{-1}BA^{-2}\rule{0ex}{4ex}&(1,0,0)
\end{array}
\>\>\>
\vline
\>\>\>
\begin{array}{lll}
\abcd{1 }{   0}{ 6}{  1    }&=A^{-1}T    \rule{0ex}{4ex}&(1,0,1)\\
\abcd{31}{ -12}{168}{ -65  }&=(ATA)^{-1}B\rule{0ex}{4ex}&(0,1,1)\\
\abcd{-29}{ 12}{-162}{ 67  }&=A^{-1}T^{-2}B  \rule{0ex}{4ex}&(1,1,0)\\
\abcd{-71}{ 30}{ -258}{ 109}&=A^{-1}BTB     \rule{0ex}{4ex}&(1,0,1)\\
\abcd{37}{-30}{132}{-107}   &=(ATB^{-1}A)^{-1}  \rule{0ex}{4ex}&(0,1,1)\\
\abcd{ 7}{-6}{6}{ -5}       &=TA^{-1}          \rule{0ex}{4ex}&(1,0,1)\\
\rule{0ex}{4ex}
\end{array}
$$

To see if a lift $\Gamma$ has level $6$, we need to check whether all the
above matrices are in $\Gamma$.
We use the algorithm described in Proposition~\ref{Prop:algorithm}.
The space
$V_{\Gamma_1(6),6}$
 spanned by the $v_i$ given in (\ref{eqn:v_i_in_alg})
is spanned by $(1,1,0)$ and $(1,0,1)$.
The matrices $A,B,C$ generate $\Gamma_1(6)$ which is a lift of
$\overline{\Gamma_1(6)}$ not containing $-1$, so the $\epsilon_i$
in (\ref{eqn:liftingcondition}) are all $0$
(as discussed in Remark~\ref{rem:notesonalg}).
Thus level $6$ congruence lifts of $\overline{\Gamma_0(6)}$
not containing $-1$
correspond to solutions of
\begin{equation}
\begin{pmatrix}
1 & 1 & 0\\
1 & 0 & 1\\
\end{pmatrix}
\begin{pmatrix}
\delta_1 \\ \delta_2\\ \delta_3
\end{pmatrix}
=
\begin{pmatrix}
0 & 0 & 0\\
\end{pmatrix}.
\end{equation}
The only possible solutions in $\F_2^3$ are
$(0,0,0)$ and $(1,1,1)$, and so
we find that
the only possible congruence level $6$ lifts of $\overline{\Gamma_1(6)}$
not containing $-1$ are
$\langle A, B, T\rangle$
and
$\langle -A, -B, -T\rangle$.

Using Magma \cite{magma}, we find that
$\overline{\Gamma(12)}$ can be generated by $97$ matrices.
However, up to parity and rearranging of the letters $A,B,T$ as above,
we find that the space $V_{\Gamma_1(6),12}$ spanned by the $v_i$ in
(\ref{eqn:v_i_in_alg})
is one dimensional, spanned by $(0,1,1)$,
corresponding for example to
$$
TBA^{2}=\abcd{169}{-36}{108}{-23}.
$$
In consequence, the groups
$\langle A,B,C\rangle,\langle A,-B,C\rangle,
\langle-A,B,-C\rangle,\langle-A,-B,-C\rangle,$
all contain $\Gamma(12)$.  The remaining lifts are noncongruence
subgroups.

The following table shows the data used to compute the dimension of
the spaces of weight $3$ cusp forms.
The last line displays whether the group is a congruence subgroup or not,
and if so, its level.
$$
\begin{array}{|lllllll|}
\hline
\text{cusp}&\text{stabilizer}&\!\!\!\!\!\!\!\!\!\!\!\!\text{width}
& \multicolumn{4}{c|}{\text{regular in
...?}}\\
&&&
\langle T,A,B\rangle&
\langle T,-A,B\rangle&
\langle -T,A,B\rangle&
\langle -T,-A,B\rangle\\
\infty & T & 1& Y & Y & n & n\\
0      & AT=(T^{-6})^{S}&6 & Y & n & n & Y\\
\frac12    & B              & 3 & Y & Y & Y & Y\\
\frac13    & A^{-1}B        & 2 &Y  &n  & Y & n\\
\multicolumn{3}{|r}{\dim S_3(\Gamma)}&0 & 1 & 1 & 1\\
\multicolumn{3}{|r}{\text{level if congruence}}& 6 & 12 & -  & -\\
\hline
\end{array}
$$

$$
\begin{array}{|lllllll|}
\hline
\text{cusp}&&& \multicolumn{4}{c|}{\text{regular in
...?}}\\
&&&
\langle T,A,-B\rangle&
\langle T,-A,-B\rangle&
\langle -T,A,-B\rangle&
\langle -T,-A,-B\rangle\\
\infty     & && Y & Y & n & n\\
0          & & & Y & n & n & Y\\
\frac12    & & & n & n & n & n\\
\frac13    & & &n & Y & n&Y\\
\multicolumn{3}{|r}{\dim S_3(\Gamma)}&1 & 1 & 2 & 1\\
\multicolumn{3}{|r}{\text{level if congruence}}&-&-& 12&6\\
\hline
\end{array}
$$

\end{ex1}

\begin{ex1}
$\Gamma_0(8)$ can be generated by the following matrices:
$$
g_1=\abcd{ 1}{ 1}{ 0 }{ 1 },\;
g_2=\abcd{ 5}{-1}{ 16}{ -3},\;
g_3=\abcd{ 5}{-2}{ 8 }{ -3}.
$$
Let $V_8=V_{\Gamma_0(8),8}$ be the vector space spanned by the $v_i$ of
(\ref{eqn:v_i_in_alg}), for $\Gamma(8)$, and let
$V_{16}=V_{\Gamma_0(8),16}$
be the space spanned by the $v_i$ for $\Gamma(16)$.
Then $V_8$ is spanned by $(1,1,1)$ and $V_{16}=\{0\}$.
So $4$ of the eight lifts not containing $-1$ are congruence subgroups
of level $8$, and the rest are congruence subgroups of level $16$.

\end{ex1}

\begin{ex1}
$\Gamma_0(16)$ is generated by the following matrices:
$$
g_1=\abcd{ 1}{ 1   }{ 0}{ 1    },\;
g_2=\abcd{ 5}{ -1 }{ 16}{ -3 },\;
g_3=\abcd{ 25}{ -9 }{ 64}{ -23 },\;
g_4=\abcd{ 9}{ -4 }{ 16}{ -7 },\;
g_5=\abcd{ 13}{ -9 }{ 16}{ -11 } .
$$
Let $V_{32}=V_{\Gamma_0(16),32}$ be the vector space spanned by the $v_i$ of
(\ref{eqn:v_i_in_alg}), for $\Gamma(32)$, and let
$V_{16}=V_{\Gamma_0(16),16}$
be the space spanned by the $v_i$ for $\Gamma(16)$.
Then
$$
V_{16}=\langle
  ( 0, 0, 0, 1, 0 ),
  ( 0, 1, 0, 0, 1 ),
  ( 1, 0, 1, 1, 0 )
\rangle
$$
and
$$V_{32}=\langle
(0, 1, 0, 0, 1 ), ( 1, 0, 1, 0, 0)
\rangle
$$
So $4$ of the $32$ lifts not containing $-1$ are congruence subgroups
of level $16$, and $4$ are congruence subgroups of level $32$.
The remaining $24$ lifts are noncongruence subgroups.

\end{ex1}

\begin{ex1}

$\Gamma_0(20)$ can be generated by the following matrices:
$$
\begin{array}{llllll}
  g_1 &=& \begin{pmatrix} 1 & 1   \\ 0 & 1    \end{pmatrix}&  g_5 &=& \begin{pmatrix} 31 & -9 \\ 100 & -29\end{pmatrix}\\
  g_2 &=& \begin{pmatrix} 13 & -2 \\ 20 & -3  \end{pmatrix}&  g_6 &=& \begin{pmatrix} 17 & -6 \\ 20 & -7  \end{pmatrix}\\
  g_3 &=& \begin{pmatrix} 31 & -7 \\ 40 & -9  \end{pmatrix}&  g_7 &=& \begin{pmatrix} 11 & -5 \\ 20 & -9  \end{pmatrix}\\
  g_4 &=& \begin{pmatrix} 29 & -8 \\ 40 & -11 \end{pmatrix}\\
\end{array}
$$
Let $V_{40}=V_{\Gamma_0(20),40}$ be the vector space spanned by the $v_i$ of
(\ref{eqn:v_i_in_alg}), for $\Gamma(40)$, and let
$V_{20}=V_{\Gamma_0(20),20}$
be the space spanned by the $v_i$ for $\Gamma(20)$.
Since $V_{40}\subset V_{20}$, we can define a space $W_{20}$
with $V_{20}=V_{40}\oplus W_{20}$.
The following table displays some of the data computed in order
to determine the number of congruence lifts of
$\overline{\Gamma_0(20)}$.
$$
\begin{array}{llll}
\text{basis of }V_{40}&
\text{basis of }W_{20}& \text{basis of }V_{20}^{\perp}
&\text{basis of }V_{40}^{\perp}\\
\begin{array}{l}
  (1, 0, 1, 0, 0, 0, 0)  \\
  (0, 1, 0, 0, 0, 1, 0) \\
  (0, 0, 0, 1, 0, 0, 0) \\
  (0, 0, 0, 0, 1, 0, 1)
\end{array}
&
\begin{array}{l}
 (1, 1, 0, 1, 1, 1, 0) \\
 (1, 0, 0, 0, 0, 0, 1)
\end{array}
&
\begin{array}{l}
(1,0,1,0,1,0,1)\\
(0,1,0,0,0,1,0)
\end{array}
&
\begin{array}{l}
(1,0,1,0,0,0,0)\\
(0,0,0,0,1,0,1)\\
(0,1,0,0,0,1,0)
\end{array}
\end{array}
$$
For example, the vector $(0,0,0,0,1,0,1)$ in $V_{40}$
corresponds to
$$\abcd{-9679}{2800}{-280}{81}=g_1^{34}g_7g_5.$$
From the above data, we see that there are $2^2=4$ congruence lifts of
level $20$ not containing $-1$, and
a further $4$ of level $40$.  The remaining $2^7-8=120$
lifts not containing $-1$ are noncongruence.
\end{ex1}

\section{Proof of Theorem \ref{Thm:gamma}}\label{Proof:gamma}

In this section, we study the principal congruence groups $\overline{\Gamma(N)}$ in $\PSL_2(\Z)$ and their lifts to $\SL_2(\Z)$. Note that for $N=1$, i.e.~$\PSL_2(\Z)$ there is only one lift, the full group $\SL_2(\Z)$, since $\pm\begin{pmatrix}0 & -1\\1 & 0 \end{pmatrix}$ has square $-1$.

We refer to $\Gamma(N)$ as the \emph{canonical lift}. For an element $\bar A\in\overline{\Gamma(N)}$ we have a unique canonical lift $A_0\in\Gamma(N)$ if $N>2$.
With respect to the canonical lift $A_0$, we define the \emph{sign} $\sigma(A)$ of any lift $A$ of $\bar A$ as
\[
 \sigma(A)=\begin{cases}
            \phantom{-}1, & \text{if } A=A_0,\\
-1, & \text{if } A=-A_0.
           \end{cases}
\]

Let $\Gamma$ denote any lift of $\overline{\Gamma(N)}$. If $\Gamma$ is congruence, then $\Gamma(2N)\subset\Gamma$ by Proposition \ref{Prop:2N}. As a consequence, for all $A\in\Gamma$ with  $\bar A\in\overline{\Gamma(2N)}$, we obtain $\sigma(A)=1$, i.e.~the lift in $\Gamma$ is fixed as $A=A_0$. Equivalently, for $A,B\in\Gamma$
\[
 \bar A \equiv \bar B \mod 2N \;\;\Longrightarrow \;\; \sigma(A)=\sigma(B).
\]
(We employ the  convention that $\bar A\equiv \bar B\mod 2N$ if and only if $\bar A\bar B^{-1}\in\overline{\Gamma(2N)}$.)

\medskip

We now specialize to the situation where $N$ is \textbf{odd}. Hence for $A, B\in\Gamma$
\[
 \bar A \equiv \bar B \mod 2N \;\;\Longleftrightarrow \;\; A \equiv B \mod 2.
\]
Thus we consider the projection from $\Gamma$ to $\SL_2(\F_2)$ (which is a group homomorphism factoring through the quotient $\overline\Gamma$):
\begin{eqnarray*}
 \Gamma & \to & \SL_2(\F_2) \\
 A & \mapsto & \tilde A
\end{eqnarray*}

The congruence property of $\Gamma$ gives the implication
\begin{eqnarray}\label{eq:2}
\tilde A=\tilde B\;\; \Longrightarrow \;\; \sigma(A)=\sigma(B).
\end{eqnarray}

There is one prominent matrix in $\Gamma$: the lift of $\bar T^N$, which we shall denote by $T_N$. Apparently its sign can be chosen freely.

\begin{lem}\label{Lem:sign}
Let $N$ be odd. If $\Gamma$ is congruence and $-1\not\in\Gamma$, then the sign of $T_N$ fixes the signs of all other elements by (\ref{eq:2}).
\end{lem}

\begin{proof}
 Let $A\in\Gamma$. We shall use the fact that $\SL_2(\F_2)\cong S_3$.
If $\tilde A$ has odd order, then $\tilde A^3=1$. 
Hence $\sigma(A)=\sigma(A^3)=1$.
If $\tilde A$ has order two, then either $\tilde A = \tilde T_N$ or $\tilde A\,\tilde T_N$ has order three. In each case $\sigma(A)=\sigma(T_N)$.
\end{proof}

\begin{cor}
 If $N>1$ is odd, then $\overline{\Gamma(N)}$ has exactly three congruence lifts.
\end{cor}

Namely the congruence lifts are $\Gamma(N), \{\pm 1\}\cdot \Gamma(N)$ and the subgroup generated by $-T^N$ and the appropriate lifts of the other generators of $\overline{\Gamma(N)}$ as determined by Lemma \ref{Lem:sign}. By construction, each of these lifts contains $\Gamma(2N)$.

\begin{cor}\label{cor:odd}
If $N>1$ is odd, then $\overline{\Gamma(N)}$ has noncongruence lifts.
\end{cor}

Since $\overline{\Gamma(N)}$ cannot be generated by a single element, the corollary follows readily. For $N=p$, we obtain explicit numbers from Frasch's work \cite{frasch}: $\overline{\Gamma(p)}$ is freely generated by $r=1+p\,(p^2-1)/12$ elements, so we have three congruence lifts and $2^r-2$ noncongruence lifts.
\medskip

We now turn to \textbf{even} level $N$. Here $\Gamma(N)/\Gamma(2N)\cong (\Z/2)^3$ by virtue of the following construction: Let $A\in\Gamma(N)$. Write $A$ uniquely as
\[
 A = 1 + N\cdot B + 2N \cdot C,\;\;\; B\in M_2(\F_2),\; C\in M_2(\Z).
\]
Here $1=\det(A)=1+N\,\mbox{tr}(B)+2N\,(\hdots)$. In consequence, $B$ has even trace. This defines a group homomorphism
\begin{eqnarray*}
\Gamma(N) & \to & M_2(\F_2)\\
 A & \mapsto & \hat A=B
\end{eqnarray*}
with kernel $\Gamma(2N)$ and image the matrices in $M_2(\F_2)$ with zero trace.
Hence it
identifies the quotient $\Gamma(N)/\Gamma(2N)$ with the image, i.e.~abstractly with $(\Z/2)^3$.

First we consider the case $N>2$. Since $-1\not\in\Gamma(N)$, there are canonical lifts. Hence for any lift $\Gamma$ of $\overline{\Gamma(N)}$, we obtain a group homomorphism
\[
 \Gamma\ni A \mapsto \bar A \mapsto A_0 \mapsto \hat A = \hat A_0.
\]

As before, if $\Gamma$ is congruence, then it contains $\Gamma(2N)$ by Proposition \ref{Prop:2N}. Hence we have
\[
 \hat A = \hat B\;\; \Rightarrow \;\; \bar A \equiv \bar B \mod 2N \;\; \Rightarrow \;\; \sigma(A) = \sigma(B).
\]
It follows that the signs of three elements with independent image in $M_2(\F_2)$ determine the congruence lift $\Gamma$.
Here two signs are given by the elements $T_N, T_N'\in\Gamma$.
The third sign is fixed by the lift of
\[
 \begin{pmatrix}
  1+N & -N\\
N & 1-N
 \end{pmatrix}.
\]

As in the case of odd $N$, the sign restrictions determine all congruence lifts of $\overline{\Gamma(N)}$:

\begin{lem}
\label{lem:even}
Let $N>2$ be even. Then $\overline{\Gamma(N)}$ has exactly 9 congruence lifts.
\end{lem}

\begin{cor}
\label{Cor:non-cong}
 If $N>2$, then $\overline{\Gamma(N)}$ has noncongruence lifts.
\end{cor}

\begin{proof}
For odd $N>1$, this result is Corollary \ref{cor:odd}.
For even $N>2$, it will follow from Lemma \ref{Lem:sub}, once we verify the claim for $N=4$.
To this end, it suffices to check that $\overline{\Gamma(4)}$ has at least four free generators.
In fact, we have seen in Example \ref{ex:4} that $\overline{\Gamma(4)}$ has five free generators.
\end{proof}


To complete the proof of Theorem \ref{Thm:gamma}, we are concerned with the remaining case $N=2$. Here the situation differs since $-1\in\Gamma(2)$. Nonetheless we can define a non-trivial canonical lift $\Gamma_0\subsetneq\Gamma(2)$ by requiring
\[
 a\equiv d\equiv 1 \mod 4
\]
for the diagonal entries of any matrix in $\Gamma_0$. Then we proceed as before with the only difference that the image of $\Gamma_0$ in $M_2(\F_2)$ consists of all matrices with zero diagonal elements. Abstractly, the image is isomorphic to $(\Z/2\Z)^2$. This agrees with the generators of $\Gamma(2)$ being $-1, T^2, (T^2)'$  by \cite{frasch}. Hence any non-trivial lift $\Gamma$ of $\overline{\Gamma(2)}$ is generated by the lifts $T_2, T_2'$. By the above consideration, for $\Gamma$ to be congruence we can choose the signs of both generators freely. Hence all five lifts are congruence.

Alternatively, one could deduce this claim from the fact that $\Gamma(2)^2 = \Gamma(4)$ as we will show in section \ref{squares}.

\section{Proof of Theorem \ref{Thm:gamma0}}\label{Proof:gamma0}

Let us first recall some well-known facts about elements of finite order of $\SL_2(\Z)$, cf.\ \cite{shimura}, Propositions 1.12 and 1.18: an element $\left( \begin{smallmatrix} a & b \\ c & d \end{smallmatrix} \right) \in\SL_2(\Z)$ different from $\pm 1$ has finite order if and only if $|a+d|<2$. An easy computation then shows that $-1$ is the only element of order $2$, whereas the elements of order $4$ are precisely the elements of form
$$
\left( \begin{array}{cc} a & b \\ c & -a \end{array} \right)
$$
where $a,b,c\in\Z$ with $a^2 + bc = -1$.
\medskip

\begin{proof}[Proof of (i)]
Write $N=2^s \cdot p_1^{s_1} \cdots p_t^{s_t}$ where $s\le 1$ and the $p_i$ are distinct primes each congruent to $1$ modulo $4$.

Thus for each $i$, the number $-1$ is a square in $\Z_{p_i}$, i.e.~the congruence $x^2 \equiv -1 \pod{p_i^r}$ is solvable in $\Z$ for every $r\in\N$. Since $s\le 1$ the same holds trivially true for the congruence $x^2 \equiv -1 \pod{2^s}$. By the Chinese Remainder Theorem we conclude that there is $a\in\Z$ such that $a^2 \equiv -1 \pod{N}$.

Consequently there exists an element $\gamma = \left( \begin{smallmatrix} a & b \\ N & -a \end{smallmatrix} \right)$ of order $4$ in $\Gamma_0(N)$. If now $\Gamma$ is a lift of $\overline{\Gamma_0(N)}$ we must have  $\gamma\in\Gamma$ or $-\gamma\in\Gamma$. Since the square of both $\gamma$ and $-\gamma$ is $-1$, we conclude that $\Gamma$ contains $-1$ and hence equals $\Gamma_0(N)$.
\end{proof}


\begin{proof}[Proof of (ii)]
Let $p$ be an odd prime, $p\equiv 3\pod{4}$ and $N=p^r$ for some $r\in\N$.
Denote by $(\cdot/p)$ the Legendre symbol modulo $p$.
We define the sign homomorphism
\begin{eqnarray*}
\sigma:\;\;
\Gamma_0(N)\; & \to & \{\pm 1\}\\
\begin{pmatrix}
a & b\\
c & d
\end{pmatrix}
& \mapsto &
\left(\frac ap\right)
\end{eqnarray*}
Then $G_1=\ker(\sigma)$ is a congruence subgroup of $\Gamma_0(N)$, since
$\Gamma_1(N)\subset G_1$.
The subgroup $G_1$ consists exactly of those matrices in $\Gamma_0(N)$ whose diagonal entries are squares modulo $p$ (or equivalently modulo $N$).

Let $\Gamma$ be a congruence lift of $\overline{\Gamma_0(N)}$. By Proposition \ref{Prop:2N} we have then $\Gamma(2N) \le \Gamma$.
Note that $\sigma$ is trivial on $\Gamma(2N)$.
Hence $\sigma$ factors through the quotient $\Gamma / \Gamma(2N)$.

In order to study this quotient, consider the homomorphism
$$
\phi : ~ \SL_2(\Z) / \Gamma(2N) \longrightarrow \SL_2(\Z/2) \times \SL_2(\Z /N)
$$
given by $A \mapsto (A \bmod{2}, A \bmod N)$;
$\phi$ is an isomorphism: it is clearly injective, and hence surjective by comparison of orders (cf.\ for instance \cite{miyake}, Theorem 4.2.5, for formulas for indices of the various subgroups in $\SL_2(\Z)$).
Under this isomorphism, $\Gamma_0(N)$ is mapped to
$$
\phi(\Gamma_0(N)) =: G = \SL_2(\Z /2) \times H \le \SL_2(\Z/2) \times \SL_2(\Z /p)
$$
where $H \cong \Gamma_0(N) / \Gamma(N)$: clearly the image is contained in $\SL_2(\Z /2) \times H$ and hence equals this group, again by order considerations:
$$
[\Gamma_0(N) : \Gamma(2N)] = 6p^{2r-1}(p-1) = 6 \cdot [\Gamma_0(N) : \Gamma(N)] = \# (\SL_2(\Z/ 2)) \cdot [\Gamma_0(N) : \Gamma(N)] .
$$
For $A\in\Gamma_0(N)$, $\sigma(A)$ is determined by the image of $A$ in $H$.
Hence we shall study the structure of $H$ to some extent.
In the present situation, $H$ is not cyclic, but close enough to it in the following sense:

\begin{lem}
\label{lem:N}
For any $A\in H$, we have $A^u=\pm 1$ for $u=p^{2r-1}(p-1)/2$.
\end{lem}

\begin{proof}
We will use the immediate fact for $A\in\Gamma_0(N)$ with diagonal entries $a, d$
that $A^t$ has diagonal entries $a^t, d^t$ modulo $N$.
Here $N=p^r$ so that $(\Z/N\Z)^*$ is cyclic of order $p^{r-1}(p-1)$.
Thus for any $A\in\Gamma_0(N)$, $A^{p^{r-1}(p-1)/2}$ has diagonal entries $\pm 1$.
But then for any $B\in\Gamma_0(N)$ with diagonal entries $\pm 1$,
we have $B^{p^r}=\pm 1$.
\end{proof}

\begin{cor}
\label{cor:N}
For any $A\in\Gamma_0(N)$, we have
\[
\sigma(A)=
\begin{cases}
1 & \text{if $A$ has odd order in $H$};\\
-1 & \text{if $A$ has even order in $H$}.
\end{cases}
\]
\end{cor}

\begin{proof}
Since $u=p^{2r-1}(p-1)/2$ is odd,
we have $\sigma(A)=\sigma(A^u)$ where $A^u=\pm 1$ in $H$ depending on the parity of the order.
The distinction of the two cases follows directly.
\end{proof}

We return to the congruence lift $\Gamma$ of $\Gamma_0(N)$.
Suppose that $-1\notin \Gamma$. Then, since $-1$ is not in the kernel of $\phi$, we have that $\phi(\Gamma)$ is a subgroup of index $2$ of $G$, and $G$ is generated by $\phi(G)$ and $\phi(-1)$. In particular $\phi(-1) = (1, -1)\not\in \phi(G)$.

\begin{prop}
\label{prop:N}
Let $\Gamma\neq G_1, \Gamma_0(N)$ be a congruence lift.
Then $\sigma$ is determined on $\Gamma$ by the relation
\[
\sigma(A)=-1 \Longleftrightarrow A \text{ has order two in } \SL_2(\Z/2).
\]
\end{prop}

\begin{proof}
Since $\Gamma\neq G_1$,
there is $B\in\Gamma$ such that $\sigma(B)=-1$.
Let $C:=B^u$.
Then $\sigma(C)=-1$, so that $C=-1$ in $H$ by Lemma \ref{lem:N} and Corollary \ref{cor:N}.
But then since $(1,-1)\not\in\phi(G)$, we deduce that $C$ has order two in $\SL_2(\Z/2)$.

First we let $A\in\Gamma$ have odd order in $\SL_2(\Z/2)$.
Then $\phi(A^{3u})=(1,1)$ since $A^{3u}=\pm 1$ in $H$, but $(1,-1)\not\in\phi(\Gamma)$.
Hence $\sigma(A)=\sigma(A^{3u})=1$.

Now let $A\in\Gamma$ have  order two in $\SL_2(\Z/2)$.
In consequence $A\cdot C$ has odd order in $\SL_2(\Z/2)$.
By the first alternative, we obtain $1=\sigma(A\cdot C)=-\sigma(A)$.
\end{proof}

We are now ready to prove the second statement from Theorem \ref{Thm:gamma0}.
Namely the order in $\SL_2(\Z/2)$ does not depend on the lift from $\PSL_2(\Z)$ to $\SL_2(\Z)$.
Hence the sign condition in Proposition \ref{prop:N} determines the lift $\Gamma\neq G_1, \Gamma_0(N)$ uniquely.
This lift $\Gamma$ is in fact congruence
as it contains $\Gamma(2N)$ by definition.
\end{proof}

\begin{proof}[Proof of (iii)]

We start by utilizing results of Rademacher to treat the prime case $N=p$:

\begin{lem}
\label{lem:number}
If $p\equiv 3\, (4)$ with $p>3$, then there are noncongruence lifts of $\overline{\Gamma_0(p)}$. More precisely, putting
$$
s := 2\left[ \frac{p}{12} \right] + 3
$$
the number of noncongruence lifts of $\overline{\Gamma_0(p)}$ is:
$$
\left\{ \begin{array}{ll} 0 & \mbox{if $p=3$} \\ 2^{s-2} - 2 & \mbox{if $p\equiv 7 \pod{12}$} \\ 2^s - 2 & \mbox{if $p\equiv 11 \pod{12}$}. \end{array} \right.
$$

\end{lem}

\begin{proof}
By \cite[pp.~146--147]{rademacher}, we know that if $p>3$ then the group $\overline{\Gamma_0(p)}$ is generated by
$$
s := 2 \left[ \frac{p}{12} \right] + 3
$$
elements and relations as follows. If $p\equiv 11\pod{12}$ there are no relations, whereas if $p\equiv 7 \pod{12}$ there are $2$ relations involving two of the generators, call them $V_1$ and $V_2$, namely
$$
V_1^3 = V_2^3 = 1 ~.
$$

By Lemma \ref{Lemma:lifts} we can determine the number of lifts when $p>3$:
in the cases $p\equiv 11\pod{12}$ and $p\equiv 7\pod{12}$, there are precisely $2^s$ resp.~$2^{s-2}$ lifts of $\overline{\Gamma_0(p)}$ not containing $-1$.
By Theorem \ref{Thm:gamma0} {\it (ii)} there are precisely two congruence lifts not containing $-1$.
Thus the formulas for the number of noncongruence lifts follow when $p>3$. One checks immediately that this number is positive in all cases.

As for the group $\overline{\Gamma_0(3)}$, by \cite{rademacher} it is generated by two elements $S$ and $V$ with the single relation $V^3=1$. Hence Lemma \ref{Lemma:lifts} implies that the group admits precisely two lifts not containing $-1$. Both lifts are congruence by Theorem \ref{Thm:gamma0} {\it (ii)}.
\end{proof}

Let $N$ be as in Theorem \ref{Thm:gamma0} {\it (iii)}.
We have seen in the section \ref{examples}, that there are noncongruence lifts of the groups $\overline\Gamma_0(6), \overline{\Gamma_0(16)}$ and $\overline{\Gamma_0(16)}$.
Using Theorem \ref{Thm:gamma0} {\it (ii)} one can prove the same for $\overline{\Gamma_0(9)}$, since this group has more than one free generator.
Thus, by the hypothesis and by Lemma \ref{lem:number}, there is a divisor $M$ of $N$ such that $\overline{\Gamma_0(M)}$ has a noncongruence lift.
Let now $\overline G$ be any of the groups $\overline{\Gamma_0(N)}$ and $\overline{\Gamma_1(N)}$.
Then $\overline G \le \overline{\Gamma_0(M)}$,
and the claim follows from Lemma \ref{Lem:sub}.
%
\end{proof}

Theorem \ref{Thm:gamma0} leaves essentially three cases unanswered:
$3, 4$ or $8$ times a product of primes congruent to $1$ modulo $4$.
An analysis of this case requires different techniques that we hope to address elsewhere.

\section{Subgroups with the same cusp forms}
\label{s:MF}

Suppose that $\Gamma_1$ and $\Gamma_2$ are subgroups of finite index in $\SL_2(\Z)$.
Let $S_k(\Gamma_i)$ denote the space  of cusp forms of weight $k$ with respect to $\Gamma_i$.
If the groups contain $-1$, then $S_k(\Gamma_1) = 0 = S_k(\Gamma_2)$ for all odd $k$. In particular this holds for infinitely many $k$.

Now consider the case where the groups do not contain $-1$.
We will show that $S_k(\Gamma_1) = S_k(\Gamma_2)$ holds for infinitely many $k$ only if the groups are projectively equivalent. More precisely:

\begin{prop}\label{prop:sameforms} (i). Suppose that $G$ and $\Gamma$ are subgroups of $\SL_2(\Z)$ of finite index and that $\Gamma$ is a subgroup of $G$.
\smallskip

Suppose that we have $\dim S_k(G) = \dim S_k(\Gamma)$ for infinitely many positive integers $k$, and that either infinitely many of these $k$ are even, or that $-1\not\in\Gamma$.
\smallskip

Then $\overline G = \overline\Gamma$.
\smallskip

\noindent (ii). Let $\Gamma_1$ and $\Gamma_2$ be subgroups of $\SL_2(\Z)$ of finite indices and suppose that $-1$ is not in $\Gamma_1 \cap \Gamma_2$.
\smallskip

Then $S_k(\Gamma_1) = S_k(\Gamma_2)$ for infinitely many $k\in\N$ if and only if $\overline\Gamma_1 = \overline\Gamma_2$.
\end{prop}
\begin{proof}[Proof of (i)] Let $\gamma$ denote the genus of the modular curve $\Gamma \backslash \mathcal{H}^{\ast}$. Here as usual, $\mathcal{H}^{\ast}$ denotes the upper halfplane with the cusps $\Q \cup \{ \infty \}$ added. Let $\mu:=[\PSL_2(\Z):\overline\Gamma]$, let $\nu_2$ and $\nu_3$ denote the number of inequivalent elliptic points of $\Gamma$ of orders $2$ and $3$, respectively, let $\nu_{\infty}$ be the number of inequivalent cusps of $\Gamma$, and let $\nu_{\infty}^+$ and $\nu_{\infty}^-$ denote the number of inequivalent regular and irregular cusps of $\Gamma$, respectively. Thus, $\nu_{\infty}=\nu_{\infty}^+ + \nu_{\infty}^-$.
\smallskip

If $k$ is even and $>2$ we have by the dimension formula, cf.\ \cite{miyake}, Theorem 2.5.2, that
$$
\dim S_k(\Gamma) = (k-1)(\gamma-1) + \frac{k-2}{2} \cdot \nu_{\infty} + \left[ \frac{k}{4} \right] \cdot \nu_2 + \left[ \frac{k}{3} \right] \cdot \nu_3
$$
and also
$$
\gamma = 1 + \frac{\mu}{12} - \frac{\nu_2}{4} - \frac{\nu_3}{3} - \frac{\nu_{\infty}}{2}
$$
by the genus formula, cf.\ \cite{shimura}, Proposition 1.40. Thus,
$$
\dim S_k(\Gamma) = \frac{k-1}{12} \cdot \mu - \frac{\nu_{\infty}}{2} + \left( \left[ \frac{k}{4} \right] - \frac{k-1}{4} \right) \cdot \nu_2 + \left( \left[ \frac{k}{3} \right] - \frac{k-1}{3} \right) \cdot \nu_3
$$
for $k$ even and $>2$.

Similarly, if $-1\not\in\Gamma$ and $k\ge 3$ is odd we have the dimension formula \cite{miyake}, Theorem 2.5.3:
$$
\dim S_k(\Gamma) = (k-1)(\gamma-1) + \frac{k-2}{2} \cdot \nu_{\infty}^+ + \frac{k-1}{2} \cdot \nu_{\infty}^- + \left[ \frac{k}{4} \right] \cdot \nu_2 + \left[ \frac{k}{3} \right] \cdot \nu_3
$$
which combines with the genus formula to:
$$
\dim S_k(\Gamma) = \frac{k-1}{12} \cdot \mu - \frac{\nu_{\infty}^+}{2} + \left( \left[ \frac{k}{4} \right] - \frac{k-1}{4} \right) \cdot \nu_2 + \left( \left[ \frac{k}{3} \right] - \frac{k-1}{3} \right) \cdot \nu_3
$$
where we used $\nu_{\infty}=\nu_{\infty}^+ + \nu_{\infty}^-$.
\smallskip

Now note that we have similar formulas for $\dim S_k(G)$, that for any $a\in\N$ the number $\left[ \frac{k}{a} \right] - \frac{k-1}{a}$ stays bounded for $k\rightarrow \infty$, and that if we have $-1\not\in\Gamma$ and $\dim S_k(G) = \dim S_k(\Gamma)$ for infinitely many odd $k$, then necessarily $S_k(G) \not= 0$ for some odd $k$ and thus $-1\not\in G$. Combining our hypothesis with an asymptotic consideration then shows that we necessarily have
$$
[\PSL_2(\Z):\overline\Gamma] = \mu = [\PSL_2(\Z):\overline G]
$$
and hence $\overline G = \overline\Gamma$ since $\overline G \le \overline\Gamma$.
\smallskip

\noindent {\it Proof of (ii).} If $\overline\Gamma_1 = \overline\Gamma_2$, then $S_k(\Gamma_1) = S_k(\Gamma_2)$ for all even $k\in\N$.
\smallskip

Conversely, suppose that $S_k(\Gamma_1) = S_k(\Gamma_2)$ for infinitely many $k\in\N$.
Consider the group $G$ generated by $\Gamma_1$ and $\Gamma_2$.
If $k$ is such that $S_k(\Gamma_1) = S_k(\Gamma_2)$,
then $S_k(G) = S_k(\Gamma_i)$ for $i=1,2$.

By hypothesis, $-1\not\in\Gamma_1$ or $-1\not\in\Gamma_2$, say $-1\not\in\Gamma_1$. Then we obtain $\overline\Gamma_1 = \overline G$ by {\it (i)}.

If $-1\not\in\Gamma_2$, then $\overline\Gamma_2=\overline G$ by {\it (i)} as claimed. If $-1\in\Gamma_2$, then $S_k(\Gamma_2) = 0$ for odd $k$. On the other hand, as $-1\not\in\Gamma_1$ we have $S_k(\Gamma_1) \not= 0$ for all sufficiently large odd $k$. Thus, our hypothesis would imply $S_k(G) = S_k(\Gamma_1) = S_k(\Gamma_2)$ for infinitely many even $k$. Hence $\overline\Gamma_2 = \overline G$ by {\it (i)}.

Thus, in any case, $\overline\Gamma_1 = \overline G = \overline\Gamma_2$, as desired.
\end{proof}

\section{Squares of congruence subgroups}\label{squares}

We note a consequence for the subgroup $\Gamma^2$ generated by squares in a given subgroup $\Gamma$.
%
From Theorem \ref{Thm:gamma}, we can deduce immediately that $\Gamma(N)^2$ is not congruence when $N>2$: this follows from Lemma \ref{Lem:sub} since any lift of $\overline\Gamma(N)$ necessarily contains $\Gamma(N)^2$ (if $\Gamma$ is a lift and if $g\in\Gamma(N)$ then $\Gamma$ contains either $g$ or $-g$, and hence in any case $g^2$). Thus, if $\Gamma(N)^2$ is congruence, then all lifts of $\overline\Gamma(N)$ are congruence.

On the other hand, it is known that $\Gamma(1)^2$ is congruence: in fact, $\Gamma(2)\subset\Gamma(1)^2$ by \cite{Newman}. (The inclusion $\Gamma(6)\subset\Gamma(1)^2$ was previously proven by J.~R.~Smith in his 1961 thesis at Michigan State, cf.~\cite{MS}.)

This leaves open the case of $\Gamma(2)^2$ which we can solve with the techniques from section \ref{s:alg}.

\begin{lem}
\label{Lem:2^2}
$\Gamma(2)^2=\Gamma(4)$.
\end{lem}

\begin{proof}
Consider the following Farey symbol for $\Gamma(4)$.
$$
\{
\xymatrix@R=0.5ex@C=0.5ex{
-\infty    \myar^{+1}      &&
\frac{-2}{1}\myar^{-2}   &&
\frac{-3}{2}\myar^{+3}       &&
\frac{-1}{1}\myar^{+3}       &&
\frac{-1}{2}\myar^{+4}      &&
\frac{0}{1}\myar^{+4}        &&
\frac{1}{2}\myar^{+5}        &&
\frac{1}{1}\myar^{+5}        &&
\frac{3}{2}\myar^{-2}        &&
\frac{2}{1}\myar^{+1}        &&
\infty
}
\}
$$
That this is a Farey symbol for $\Gamma(4)$ follows from the fact
that all the corresponding generating matrices are in
$\Gamma(4)$, and this Farey symbol corresponds to a group with index
$24$ in $\PSL_2(\Z)$, which is equal to $[\PSL_2(\Z):\overline{\Gamma(4)}]$
using the well known formula.

We can write the generating matrices in terms of squares of elements of
$\Gamma(2)$ thus:
$$
\renewcommand\arraystretch{1.5}
\begin{array}{lll}
\abcd{1}{4}{0}{1}&=&\abcd{1}{2}{0}{1}^2\\
\abcd{-7}{12}{4}{-7}&=&(\abcd{1}{2}{0}{1}^{2})^{-1}\abcd{5}{-8}{2}{-3}^2\\
\abcd{5}{4}{-4}{-3}&=&\abcd{-3}{-2}{2}{1}^2\\
\abcd{1}{0}{4}{1}&=&\abcd{1}{0}{2}{1}^2\\
\abcd{5}{-4}{4}{-3}&=&\abcd{-3}{2}{-2}{1}^2
\end{array}
$$
This demonstrates that $\Gamma(4)\subseteq\Gamma(2)^2$. On the other hand, we can easily check that the square of any element
of $\Gamma(2)$ must be in $\Gamma(4)$ (since if $A-I=2B$ for some matrices $A, B$ with integer entries, then $A^2 - I=(2B+I)^2-I=4B(B+I)$), so $\Gamma(2)^2\subseteq\Gamma(4)$. Hence the two groups are equal.
\end{proof}

In consequence, we can identify the principal congruence subgroups whose squa\-res are congruence again:

\begin{prop}\label{Cor:square}
$\Gamma(N)^2$ is congruence if and only if $N\leq 2$.
\end{prop}

Along the same lines,
we obtain from Theorem \ref{Thm:gamma0} {\it (iii)} with Lemma \ref{Lem:sub}:

\begin{cor}
If $N$ is divisible by  $6, 9, 16, 20$ or by a prime $p>3$ congruent to $3$ modulo $4$,
then $\Gamma_0(N)^2$ and $\Gamma_1(N)^2$ are not congruence.
\end{cor}

For future use in the general case, we also note the following lemma which slightly simplifies the congruence question for squares:

\begin{lem}
If $\Gamma$ has general level $N$ and $\Gamma^2$ is congruence, then $\Gamma^2$ has level $2N$.
\end{lem}

\begin{proof}
Arguing with the cusps as before, we see that $\Gamma^2$ has general level $2N$. Hence if $\Gamma^2$ is congruence, then the level divides $4N$ by Proposition \ref{Prop:2N}.

Looking at the proof of Proposition \ref{Prop:2N}, the main difference when working with $\Gamma^2$ is that $T^{2N}\in\Gamma^2$ due to the assumption that $\Gamma$ has general level $N$. Hence Wohlfahrt's original argument goes through without doubling the level twice, since this was only necessary to ensure that $T^{2N}$ is contained in the group. Thus $\Gamma^2$, if it is congruence, has level $2N$.
\end{proof}

\begin{rem}
With the generators in the proof of Lemma \ref{Lem:2^2}, we can also give a more precise description of the lifts of $\overline{\Gamma(4)}$:
Any congruence lift of $\overline{\Gamma(4)}$ must have level either $4$ or $8$ by Proposition \ref{Prop:level}. With $V_4=V_{\Gamma(4),4}$ and $V_8=V_{\Gamma(4),8}$, as described in Remark~\ref{rem:notesonalg}, in terms of the above set of generators, we have
$$V_4^\perp=0,
\text{ and }
V_8=\langle
(0,0,1,0,1),
(1,1,0,1,0)
\rangle
$$
So there is a unique level $4$ congruence lift of $\overline{\Gamma(4)}$ not containing $-1$, namely $\Gamma(4)$, and
a further $7$ level $8$ lifts, in accordance with Lemma \ref{lem:even}; the remaining $24$ lifts are noncongruence.
\end{rem}

%








\section{Elliptic modular surfaces}\label{ss:EMS}

We conclude this paper with some geometric considerations.
Our motivation stems from modular forms.
We have seen in section \ref{s:MF} that lifts can only be distinguished on modular forms of odd weight.
For weight 3, there is an instructive relation to holomorphic 2-forms on certain complex algebraic surfaces that we shall briefly recall.

Let $\Gamma$ denote a  finite index subgroup $\Gamma$ of $\SL_2(\Z)$ not containing $-1$.
A construction by Shioda associates $\Gamma$ with  an elliptic modular surface $S(\Gamma)$ over the modular curve $X(\Gamma)$ (uniquely up to $\C$-isomorphism) \cite{ShEMS}.
For congruence subgroups, this construction simply exhibits the universal modular curve.

Throughout this paper, we have been considering different lifts.
In particular, the modular curve $X(\Gamma)$ is the same for all lifts, as it only depends on the image $\bar\Gamma$ of $\Gamma$ in $\PSL_2(\Z)$.
For shortness, we write $X=X(\Gamma), S=S(\Gamma)$.

Of course, the $j$-map
\[
 j:\;\; X \to \PP^1
\]
is also independent of the lift.
On the other hand, $S$ is determined by $j$ only up to quadratic twisting.
This notion refers to non-isomorphic elliptic curves which become isomorphic over a quadratic extension of the ground field.
In terms of an extended Weierstrass form
\[
y^2 = x^3 + Ax^2+Bx+C,\;\;\; A, B, C\in k,
\]
quadratic twists are in correspondence with squarefree elements $d\in k$.
The twist by $d$ is exhibited by the Weierstrass form
\[
y^2 = x^3 + dA^2x+d^2Bx+d^3C
\]
or equivalently
\begin{eqnarray}
\label{eq:twist}
dy^2 = x^3 + Ax^2+Bx+C.
\end{eqnarray}
The impact of a quadratic twist on the singular fibers is well understood thanks to Tate's algorithm \cite{Tate}.
The fiber type stays the same if $d$ vanishes to even order at the corresponding cusp on $X$.
In case of odd vanishing order, the fiber types change in a canonical way unless
the characteristic is three and there is wild ramification (see e.g.~\cite{MP}).
For instance, fibers of type $I_n$ are always interchanged with $I_n^*$ in Kodaira's notation.

In the case of elliptic modular surfaces, we can use some crucial extra information. Namely, any (complex) elliptic modular surface is extremal: it has maximal Picard number $\rho(S)=h^{1,1}(S)$, but finite group of sections.
By a theorem of Nori \cite{Nori}, these conditions (with $j$ non-constant) imply the absence of singular fibers of type $II, III, IV, I_0^*$.

In consequence, on elliptic modular surfaces of different lifts of $\bar\Gamma$, the quadratic twisting  can only occur at two kinds of points on $X(\Gamma)$:
\begin{itemize}
\item
cusps underneath
singular fibers of type $I_n, I_n^* (n>0)$.
Note that the local monodromy at these fibers is
\[
 T^n\;\;\; \text{ resp. } -T^n.
\]
\item
points such that $d$ vanishes to even order.
\end{itemize}
Conversely, any such twist is again modular by \cite{Nori}.
Since it has the same modular curve $X$ by construction, the projective images of the associated subgroups have to coincide.
Since at twisting points of the second kind the fiber type does not change,
the problem of detecting isomorphic elliptic surfaces is non-trivial;
Stiller gave an example of this phenomenon in \cite{Stiller} related to the commutator subgroup of $SL_2(\Z)$.
%
%
%
%
%
There are no points of the second kind if the modular curve $X(\Gamma)$ is rational.
The following two examples  go back to section \ref{s:ex}.
Equations can be found in \cite{MP} for instance, but the cusps are normalised in a different manner than in section \ref{s:ex}.

\subsection{Example: $\Gamma_1(4)$}
\label{ss:4}
The elliptic modular surface for $\Gamma_1(4)$ is rational over $\Q$ with singular fibers of types $I_1, I_4, I_1^*$. Two of the twists are rational as well, while the twist with three non-reduced fibers is K3.
In terms of Example \ref{ex:4}, this modular surface corresponds to the congruence lift $G_3$.
Recall that $S_3(G_3)$ has dimension one.
To find the normalised cusp form, we use that its square is in $S_6(G_3)=S_6(\Gamma_1(4))$.
The latter space is also one-dimensional and generated by $\eta(2\tau)^{12}$.
Hence $S_3(G_3)$ is generated by $\eta(2\tau)^6$.

This K3 surface is also modular in another sense:
its zeta function contains a factor corresponding to (a twist of) $\eta(4\tau)^6$,  the unique weight 3 newform of level $16$.
We emphasise that the two cusp forms $\eta(2\tau)^6$ and $\eta(4\tau)^6$ have Fourier expansions in terms of different uniformisers, but with the same Fourier coefficients.

The modularity of the zeta function follows from a more general result by Livn\'e \cite{Livne}.
In particular, this result applies to all singular K3 surfaces over $\Q$, i.e.~with maximal Picard number $\rho=20$.
We will see more instances of this modularity in the next example, but first we shall point out a general fact in this context.

\medskip

The above arithmetic relations with modular forms should not be seen as a surprise.
In fact, Shioda showed in \cite{ShEMS} that there is an analytic isomorphism
\begin{eqnarray}\label{eq:3}
 S_3(\Gamma) \cong H^{2,0}(S(\Gamma))
\end{eqnarray}
between cusp forms of weight 3 with respect to~$\Gamma$ and holomorphic 2-forms on $S(\Gamma)$.
In the congruence case, this isomorphism admits an arithmetic interpretation.
Namely, $S(\Gamma)$ has a model over $\Q$.
By Deligne \cite{Deligne}, the Galois representation on the cohomology of $S(\Gamma)$ splits off two-dimensional subrepresentations.
Here the traces of Frobenius correspond to eigenvalues of the Hecke operators on $S_3(\Gamma)$
as we have seen in the example above.
In the noncongruence case, there is a conjectural congruence relation due to Atkin--Swinnerton-Dyer (cf.~\cite{Long} and the references therein).

\subsection{Example: $\Gamma_1(6)$}
\label{Ex:6}

We conclude this paper with a detailed analysis of
the elliptic modular surface $S$ for $\Gamma_1(6)$.
We study the lifts of $\Gamma_1(6)$ and compute the corresponding cusp forms of weight $3$.
Then we compare with the zeta function of the corresponding twists of $S$.

A model of $S$ over $\Q$ can be given as follows:
\[
S:\;\;\; y^2 = x (x^2+(t^2-6t-3) x +16t).
\]
Either six-torsion section $(4t, \pm 4t(t-1))$ generates the Mordell-Weil group of $S$:
\[
\mbox{MW}(S)=\{(4t,\pm 4t(t-1)), (4,\pm 4(t-1)), (0,0)\}.
\]
The j-invariant of $S$ is
$$j=\frac{(t - 3)^3(t^3 - 9t^2 + 3t - 3)^3}
{(t - 9)(t - 1)^3t^2}.$$
One can describe the elliptic parameter $t$ by
$$t=9\eta(2\tau)^4\eta(3\tau)^8\eta(\tau)^{-8}\eta(6\tau)^{-4}.$$
The connection with the set-up in Example \ref{ex1:6} is summarised in the following table:
$$
\begin{array}{|c|cccc|}
\hline
\text{cusp }c & \infty & 1/3    & 1/2 &0     \\
\text{width } & 1      & 2      & 3   &6     \\
t(c)          & 9      & 0      & 1   &\infty\\
\text{stabilizer}&T    &AB^{-1} & B   &AT    \\
\hline
\end{array}
$$
Thus $S$ is a  rational elliptic surface with singular fibers of types $I_1, I_2, I_3, I_6$.
As explained, the different lifts of $\Gamma_1(6)$ from $\PSL_2(\Z)$ to $\SL_2(\Z)$ correspond
to twists of $S$ by certain $d\in\Q(t)$ as in \eqref{eq:twist}.
We collect this information in the next table.
Six of the twists are K3 surfaces as indicated.
The zeta function of each of them again contains a factor corresponding to some newform of weight three by \cite{Livne}.
The precise newform depends on the model of the K3 surface.
By the classification in \cite{S-CM}, very few Fourier coefficients suffice to determine the newform (see \cite[Rem.~2.]{S-CM}).
In the present situation, we achieved this by point counting at the first few good primes through Lefschetz' fixed point formula.

We express the resulting newforms in terms of some particular newforms (which are minimal in the sense of \cite{S-CM}) involving twists by the Legendre symbols $\chi_{u}=(u/p)$.
Let
$$f_8=\eta(\tau)^2\eta(2\tau)\eta(4\tau)\eta(8\tau)^2,
\;\;
f_{12}=(\eta(2\tau)\eta(6\tau))^3$$
and let $f_{24}$ denote the newform of level $24$ from \cite[Table 1]{S-CM}
(given by a Hecke character for $\Q(\sqrt{-6})$ of $\infty$-type $2$ and conductor $(2\sqrt{-6})$).

The full quadratic twist of $S$ is not K3,  but has $h^{2,0}=2$ and Euler number $e=36$.
In terms of the Enriques-Kodaira classification of algebraic surfaces, it is honestly elliptic (Kodaira dimension one).
We
compute the decomposition of its zeta function in terms of cusp forms $g_1, g_2$ below (see \eqref{eq}).
$$
\begin{array}{lllllll}
\hline
d          & \text{fibres}  &\text{lift of } \Gamma_1(6) & \text{level}&\dim S_3&\text{zeta}&\text{surface}\\
\hline
1        & 1^{\ },2^{\ },3^{\ },6^{\ }&\langle\;\;T,\;\;A,\;\;B\rangle&6 &0&-  &\text{rational}\\
t-9        & 1^{*},2^{\ },3^{\ },6^{*}&\langle -T,\;\;A,\;\;B\rangle  &- &1&f_{24}\otimes\chi_3 &\text{K3}\\
t(t-9)     & 1^{*},2^{*},3^{\ },6^{\ }&\langle -T, -A,\;\;B\rangle    &- &1&f_8\otimes\chi_3    &\text{K3}\\
(t-1)(t-9) & 1^{*},2^{\ },3^{*},6^{\ }&\langle -T, -A, -B\rangle      &6 &1&f_{12} &\text{K3}
\\
t(t-1)(t-9)& 1^{*},2^{*},3^{*},6^{*}  &\langle -T,\;\;A, -B\rangle    &12&2 & g_1, g_2     &e=36\\
t          & 1^{\ },2^{*},3^{\ },6^{*}&\langle\;\;T, -A,\;\;B\rangle  &12&1&f_{12}&\text{K3}\\
(t-1)      & 1^{\ },2^{\ },3^{*},6^{*}&\langle\;\;T, -A, -B\rangle    &- &1&f_8\otimes\chi_{-1}   &\text{K3}\\
t(t-1)     & 1^{\ },2^{*},3^{*},6^{\ }&\langle\;\;T,\;\;A, -B\rangle  &- &1&f_{24}\otimes\chi_2&\text{K3}
\\
\hline
\end{array}
$$

In order to compute the cusp forms in $S_3(\Gamma)$ for the lifts $\Gamma$, we shall again argue with extracting square roots, this time from $S_6(\Gamma)=S_6(\Gamma_1(6))$.
In contrast to the case in \ref{ss:4} where the dimensions matched exactly, we have to take into account the following subtlety for square roots of $f\in S_6(\Gamma_1(6))$:
For $\sqrt{f}$ to define a holomorphic function with respect to the lift $\Gamma$,
we require $f$ to vanish quadratically at all regular cusps of $\Gamma$ and at all zeroes in the upper half plane.
In the present situation, this criterion suffices to find all cusp forms.

We shall use the following modular forms of weight one for $\Gamma_1(6)$:
$$
\renewcommand\arraystretch{1.5}
\begin{array}{|lc|c|c|c|c|}
\hline
\multicolumn{2}{|r|}{\text{cusps (and widths)}} & \infty (1)& \frac13 (2) & \frac12 (3)
& 0 (6)\\
\cline{3-6}
\text{weight one forms for } \Gamma_1(6)
&&\multicolumn{4}{|c|}{\text{order of vanishing}}\\
\hline
\multicolumn{1}{|l}{a=\frac{\eta(z)\eta(6z)^6}{\eta(2z)^2\eta(3z)^3}=q - q^2 + q^3 + q^4 + \cdots }&& 1 & 0 & 0 & 0\\
\multicolumn{1}{|l}{b=\frac{\eta(2z)\eta(3z)^6}{\eta(z)^2\eta(6z)^3}=1 + 2q + 4q^2 + 2q^3 + \cdots}& & 0 & 1 & 0 & 0\\
\multicolumn{1}{|l}{c=\frac{\eta(3z)\eta(2z)^6}{\eta(6z)^2\eta(z)^3}=1 + 3q + 3q^2 + 3q^3 + \cdots}& & 0 & 0 & 1 & 0\\
\multicolumn{1}{|l}{d=\frac{\eta(6z)\eta(z)^6}{\eta(3z)^2\eta(2z)^3}=1 - 6q + 12q^2 - 6q^3\cdots}  & & 0 & 0 & 0 & 1\\
\hline
\end{array}
$$
Note that there are relations $c=a+b, d=b-8a$.
It follows that any product of six of them defines a weight $6$ form for $\Gamma_1(6)$ -- and cusp forms include the product $g=abcd=(\eta(\tau)\eta(2\tau)\eta(3\tau)\eta(6\tau))^2$ which in fact generates $S_4(\Gamma_1(6))$.
Looking at dimensions, it follows that these cusp forms generate $S_6(\Gamma_1(6))$.
A basis could be chosen as $ga^2, gab, gb^2$.

Consider the six lifts that geometrically correspond to K3 surfaces.
For each lift $\Gamma$,
we will exhibit a cusp form of weight $6$ for $\Gamma_1(6)$
whose square root gives a cusp form of weight $3$ for $\Gamma$.
As we have seen in Example \ref{ex1:6}, each lift has two regular cusps where extracting a square root requires quadratic vanishing of the weight $6$ cusp form $f$.
At the two irregular cusps, simple vanishing of $f$ suffices.
By inspection, these conditions determine the following normalised cusp forms of weight $6$ for $\Gamma_1(6)$ whose square roots yield cusp forms of weight $3$ for the six lifts:
$$
\renewcommand\arraystretch{1.5}
\begin{array}{|lc|c|c|c|c|l|}
\hline
\multicolumn{2}{|r|}{\text{cusps (and widths)}} & \infty (1)& \frac13 (2) & \frac12 (3)
& 0 (6) &\\
\cline{3-6}
\text{ weight six forms for } \Gamma_1(6)
&&\multicolumn{4}{c|}{\parbox{1.2in}{order of vanishing}} & \text{lift}\\
\hline
abc^2d^2 = \eta(z)^5\eta(2z)^5\eta(3z)\eta(6z)
&  & 1 & 1 & 2 & 2&
\langle -T,-A,B\rangle\\
a^2b^2cd =
\eta(z)\eta(2z)\eta(3z)^5\eta(6z)^5
&  & 2 & 2 & 1 & 1 &
\langle T,-A,-B\rangle\\
a^2bc^2d=\eta(2z)^6\eta(6z)^6
&  & 2 & 1 & 2 & 1 &\langle T,-A,B\rangle \\
ab^2cd^2=\eta(z)^6\eta(3z)^6
&  & 1 & 2 & 1 &2 &
\langle -T,-A,-B\rangle\\
a^2bcd^2 =\eta(z)^9\eta(2z)^{-3}\eta(3z)^{-3}\eta(6z)^9
&  & 2 & 1 & 1 &2 &
\langle T,A,-B\rangle\\
ab^2c^2d =\eta(z)^{-3}\eta(2z)^9\eta(3z)^9\eta(6z)^{-3}
&  & 1 & 2 & 2 &1 &
\langle -T,A,B\rangle\\
\hline
\end{array}
$$

In the congruence cases, the zeta function of the elliptic modular surface thus indeed agrees with that of the corresponding generator of $S_3(\Gamma)$.
In the non-congruence cases, one easily checks the Atkin--Swinnerton-Dyer relations for the square root of the corresponding weight six form generating $S_3(\Gamma)$ (which does not have integral Fourier coefficients) and the cusp form associated to $S$ (which is the indicated twist of $f_8$ or $f_{24}$).

Finally we come to the lift $\Gamma=\langle -T,\;\;A, -B\rangle$ where all cusps are irregular.
Here $S_3(\Gamma)$ is two-dimensional, and any product of two cusp forms gives a form in $S_6(\Gamma_1(6))$.
\begin{lem}\label{lem:last}
For the lift $\Gamma=\langle -T,\;\;A, -B\rangle$, we have with $\sqrt{g}=\eta(\tau)\eta(2\tau)\eta(3\tau)\eta(6\tau)$
\[
S_3(\Gamma) = \langle \sqrt ga, \sqrt gb\rangle.
\]
\end{lem}

\begin{proof}
Fix a basis $h_1, h_2$ of $S_3(\Gamma)$. A priori, we only have $h_1,h_2\in\C[[q^{1/2}]]$, but since their product is in $S_6(\Gamma_1(6))$ and thus in $\C[[q]]$, we actually deduce that either $h_1,h_2\in q^{1/2} \C[[q]]$ or $h_1,h_2\in q\C[[q]]$. Under the latter alternative we would have a cusp form $h\in q^2\C[[q]]\cap S_3(\Gamma)$. Then $h^2$, which is in $S_6(\Gamma_1(6))$, would vanish to order $4$ at $\infty$. In comparison, the cusp form in $S_6(\Gamma_1(6))$ with the highest vanishing order at $\infty$ is $ga^2$, but the vanishing order is only three, contradiction.

Under the first alternative, there is a cusp form $h\in q^{3/2}\C[[q]]\cap S_3(\Gamma)$. Along the same lines as above, one deduces that $h^2=\lambda ga^2$ for some $\lambda\in\C$. This shows that $\sqrt ga\in S_3(\Gamma)$. We proceed by considering the product $\sqrt gah_1\in S_6(\Gamma_1(6))$. In terms of our basis, we write this product as $g\gamma$ for some quadratic form $\gamma(a,b)$. By assumption, we can also write $h_1=\sqrt{g\delta}$ for some quadratic form $\delta(a,b)$. This gives an equality of modular forms $a\sqrt\delta=\gamma\in M_2(\Gamma_1(6))$. Upon squaring, we obtain the relation $a^2\delta=\gamma^2$ in $M_4(\Gamma_1(6))$. Since $M_4(\Gamma_1(6))$ has basis $a^4, a^3b,a^2b^2, ab^3, b^4$, spelling out the last relation shows that the quadratic form $\delta$ is in fact the square of a linear form in  $a,b$. As the same argument applies to $h_2$, the lemma follows.
\end{proof}


Lemma \ref{lem:last} exhibits a basis of the vector space $S_3(\Gamma)$.
Recall that the lift $\Gamma$ is congruence.
Hence the zeta function of the elliptic modular surface $S=S(\Gamma)$ encodes the Fourier coefficients of another basis $g_1, g_2$ of $S_3(\Gamma)$ which consists of Hecke eigenforms.
We compute the relevant degree $4$ factors $L_p(T)$ of the zeta function of $S$ at the first few good primes with Lefschetz' fixed point formula (over $\F_p$ and $\F_{p^2}$).
We then read off the resulting Fourier coefficients $a_p$ of $g_1, g_2$ (up to conjugacy) from the real quadratic factors of $L_p(T)$:
$$
\begin{array}{|ccc|}
\hline
p & L_p(T) & a_p\\
\hline
5 & T^4+18T^2+5^4 & \pm4\sqrt{2}\\
7 & (T^2+6T+49)^2 & -6\\
11 & T^4+210T^2+11^4 &  \pm 4\sqrt{2}\\
13 & (T^2-20T+13^2)^2 & 20\\
17 & T^4+66T^2+17^4 & \pm 16\sqrt{2}\\
\hline
\end{array}
$$
With some linear algebra on the vector space $S_3(\Gamma)$, one then verifies that (for the right sign choices between the Fourier coefficients of $g_1$ and $g_2$)
\begin{eqnarray}
\label{eq}
(g_1+g_2)/2=\sqrt{g}b,\;\;\; (g_1-g_2)/(2\sqrt{2}) = \sqrt{g}a.
\end{eqnarray}

\subsection*{Acknowledgement}
We are indebted to Andreas Schweizer for pointing out a mistake in an earlier version of this paper.
Our thanks go to the referee for several comments that helped improve the paper.

\end{document}